\documentclass{article}
\usepackage{amssymb}
\usepackage[normalem]{ulem}
\usepackage{amsfonts}
\usepackage{colortbl}
\usepackage{color}

\input{tcilatex}
\begin{document}

\begin{center}
\bigskip

\bigskip

{\Large Strategy-proofness on the Non-Paretian Subdomain}\bigskip 

Donald E. Campbell $\cdot $ Jerry S. Kelly
\end{center}

\bigskip

\qquad \qquad {\Large Abstract}\medskip 

\qquad \qquad Let $g$ be a strategy-proof rule on the domain $NP$ of profiles

\qquad \qquad where no alternative Pareto-dominates any other. Then we

\qquad \qquad establish a result with a Gibbard-Satterthwaite flavor: $g$ is 

\qquad \qquad dictatorial if its range contains at least three alternatives.
\bigskip 

\bigskip

\textbf{1. Introduction.}

\textbf{2. Framework.}

\textbf{3. Paths in NP(n, m) and an equivalence theorem.}

\textbf{4. Double induction basis: 3 individuals and 3 alternatives.}

\textbf{5. Induction on n (for m = 3).}

\textbf{6. Induction on m.}

\textbf{7. Conclusion.}

\textbf{8. References}

\bigskip

\section{Introduction\protect\medskip}

This paper shows that dictatorial rules are the only ones satisfying
strategy-proofness on the domain $NP$, the set of all profiles of strong
preferences at which no alternative Pareto-dominates any other. All existing
proofs of Gibbard-Satterthwaite employ profiles that are not in $NP$.%
\footnote{%
Of the many proofs of the Gibbard-Satterthwaite theorem on the domain of all
profiles of strong orders (e.g., Gibbard (1973), Satterthwaite (1975),
Schmeidler and Sonnenschein (1998), Barber\`{a} (1980, 1983a, 1983b, 2001),
Svensson (1999), Beno\^{\i}t (2000), Reny (2001), Sen (2001), Larsen and
Svensson (2006)), every single one makes critical use of profiles in which
some alternatives Pareto-dominate others. For example, many invoke Arrow's
theorem after first using a choice function to construct a social ordering
by seeing which of two alternatives is chosen when they are both brought to
the top of everyone's ordering, and so Pareto dominate everything else.
Aswal, Chatterji, and Sen (2003) characterize the family of domains on which
every strategy-proof social choice function with full range is dictatorial,
but only for domains that are product sets.}\medskip

There are two main reasons for studying strategy-proofness on $NP$. The
first is that profiles $u$ not in $NP$ may be quite unlikely. The society
may be large and the set of alternatives small. Then for no pair of
alternatives is a unanimous agreement likely.\medskip

Alternatively, in a two-stage group decision process the probability of
Pareto domination would be zero if the first stage narrows the set of
alternatives by rejecting, among others, alternatives that are
Pareto-dominated. The domain NP of surviving alternatives comprise the
feasible set, from which the second stage group chooses. \medskip

The second reason is that we might start with a (manipulable) rule $g$
defined on all profiles of strong preferences but where the restriction $%
g|_{NP}$ of $g$ to $NP$ is strategy-proof. Then the conclusion that $g|_{NP}$
is dictatorial will help determine the behavior of $g$. This paper stems
from a companion paper (Campbell and Kelly, 2014), which treated universally
beneficial manipulation ($UBM$) rules - defined by the property that if
anyone manipulates, everyone gains and no one is hurt. Such $UBM$ rules must
be strategy-proof on $NP$. Completing our classification analysis in that
earlier paper requires proving a Gibbard-Satterthwaite result on $NP$%
.\medskip

\medskip

This domain $NP$ is very large, even as a fraction of all profiles of strong
preferences, and so our result may not seem unexpected, but the difficulty
of our proof \textit{is} surprising. Of course, whenever we take up a proper
subdomain, we must reconsider strategy-proofness, as we will have excluded
manipulations that would have been possible on the larger domain. But
strategy-proofness is more difficult to analyze here because $NP$ is not a
product set.\medskip

It should also be noted that just because $NP$ is large is no guarantee of
dictatorship for full-range strategy-proof rules. There are supersets of $NP$
for which there do exist non-dictatorial strategy-proof rules: Suppose the
number of individuals, $n$, is odd, the set of alternatives is $\{x,y,z\}$
and the domain is $D$, the union of $NP$ and the collection of all the
profiles for which $z$ is the top-ranking alternative for every individual.
The rule on $D$ that selects the majority winner between $x$ and $y$ at
every profile in $NP$ and selects $z$ elsewhere in $D$ is non-dictatorial
and has full range. It is strategy-proof because, for example, if $x$ is
selected then someone ranks $x$ above $z$ and thus no individual who prefers 
$z$ to $x$ can unilaterally cause $z$ to be selected.\medskip

\qquad After introducing terminology and notation in Section 2, we prove
some intermediate results in Section 3. Then our main result is proved by
induction. The basis step, for three individuals and three alternatives is
established in Section 4. Induction on the number of individuals takes place
in Section 5. Then induction on the number of alternativess takes place in
Section 6.\bigskip

\section{Framework.\protect\medskip}

\qquad For given $m$, $n\geq 3$, we consider a finite set $X$ of
alternatives where $\mid X\mid =m$ and finite set $N=\{1,2,...,n\}$ of
individuals. A (strong) \emph{ordering} on $X$ is a complete, asymmetric,
transitive relation on $X$ and the set of all such orderings is $L(X)$. For $%
R\in L(X)$ and $Y\subset X$ let $R\mid Y$ denote the relation $R\cap Y\times
Y$ on $Y$, the restriction of $R$ to $Y$. If $R$ is a member of $L(X)$ we
let $R^{-1}$ denote the \emph{inverse} of R: That is, $(x,y)\in R^{-1}$ if
and only if $(y,x)\in R$. \medskip

\qquad A \emph{profile} $p$ is a map from $N$ to $L(X)$, where $%
p=(p(1),p(2),...,p(n))$ and we write $x\succ _{p(h)}y$ if individual $h$
strongly prefers $x$ to $y$ at profile $p$. The set of all profiles is $%
L(X)^{N}$. A \emph{domain} $\wp $ is subset of $L(X)^{N}$. For each subset $%
Y $ of $X$ and each profile $p$ in $\wp $ \ \ \ let $p\mid Y$ denote the
restriction of profile $p\in \wp $ to Y. That is, $p\mid Y$ represents the
function $q\mid L(Y)^{N}$ satisfying $q(i)=p(i)\mid Y$ for all $i\in N$. A 
\emph{social choice rule} on $\wp $ is a function $g:\wp \rightarrow X$,
where $\wp $ is a nonempty subset of $L(X)^{N}$. Rule $g$ is \textit{%
dictatorial} if there exists an individual $i$ such that for each profile $p$
in $\wp $, alternative $g(p)$ is the top-ranked element in $p(i)$ restricted
to Range($g$). \ A rule $g$ is \textit{full-range} if Range($g$) $=X$.

\qquad In this paper, we consider social choice rules on the \emph{%
Non-Paretian domain}, $NP(n,m)$, the set of all profiles $p\in L(X)^{N}$
such that no alternative is Pareto-dominates any other. 
\begin{eqnarray*}
NP(n,m) &=&\{p\in L(X)^{N}:\text{for every }x,y\in X,\text{if }x\neq y\text{%
, then we have} \\
x &\succ &_{p(i)}y\text{ for some }i\in N\text{ and }y\succ _{p(j)}x\text{
for some }j\in N\}
\end{eqnarray*}

It is very important to note that $NP(n,m)$ is not a Cartesian product set.
This will greatly complicate our analysis.\medskip

\qquad Two profiles $p$ and $q$ are $h$-\emph{variants}, where $h\in N$, if $%
q(i)=p(i)$ for all $i\neq h$. Individual $h$ can \emph{manipulate} the
social choice rule $g:\wp \rightarrow X$ at $p$ via $p^{\ast }$ if $p$ and $%
p^{\ast }$ belong to $\wp $, $p$ and $p^{\ast }$ are $h$-variants, and $%
g(p^{\ast })\succ _{p(h)}g(p)$. And $g$ is \emph{strategy-proof} if no one
can manipulate $g$ at any profile.\medskip

\qquad We will repeatedly use the following simple consequence of
strategy-proofness without explicitly alluding to it:\medskip

\qquad If we switch adjacent alternatives $a$ and $b$ in some individual $i$%
's preference ordering, without changing anyone else's preferences, then the
selected alternative will not change unless $a$ is selected before the
switch and $b$ is selected after and individual $i$ preferred $a$ to $b$
initially, or $b$ is selected before and $a$ is selected after and $i$
preferred $b$ to $a$ initially.\medskip

\qquad Rule $g$ on $L(X)^{N}$ satisfies $universally$ $beneficial$ $%
manipulation$ (UBM) if for every profile $u$ and individual $h$ such that
there exists an $h$-variant profile $u^{\ast }$ with $g(u^{\ast })\succ
_{u(h)}g(u)$, we have $g(u^{\ast })\succ _{u(j)}g(u)$ for every individual $%
j $. Our primary motivation for studying strategy-proofness on $NP(n,m)$ is
that if $g$ satisfies UBM on $L(X)^{N}$, then the restriction of $g$ to $%
NP(n,m)$ is strategy-proof on that subdomain. The classification theorems
(Campbell and Kelly, 2014) characterizing all UBM rules rely on the claim
that if $g:NP(n,m)\rightarrow X$ is strategy-proof and has a range of at
least three alternatives, then $g$ is dictatorial. We will exclude $n=2$,
for in that case every rule is strategy-proof on $NP(n,m)$ for the trivial
reason that no two distinct profiles in $NP(n,m)$ are $i$-variants of one
another.\bigskip

\section{Paths in NP(n, m) and an equivalence theorem.\protect\medskip}

\qquad From a global perspective, our proof will use two inductions. First,
we show strategy-proofness plus full range implies dictatorship for $m=3$
and $n=3$. We show in Section 5 that the result for $m=3$ and $n$ implies
the result for $m=3$ and $n+1$. Finally, Section 7 contains the proof that
the result for $m$ and $n$ implies the result for $m+1$ and $n$.\medskip

\qquad Those induction steps will require some preliminary theorems. One
goal of this section is to show that in our proof we can restrict attention
to rules that satisfy a full range assumption. \ Before that, we need to
show that on $NP(n,m)$, any strategy proof rule $g$ with range $S$ has the
property that $u\mid S=u^{\ast }\mid S$ implies $g(u)=g(u^{\ast })$. On $%
L(X)^{N}$, this is easy, just use a standard sequence argument in the manner
of Gibbard (1973) and Satterthwaite (1975). But here, because $NP(n,m)$ is
not a Cartesian product, those standard sequences can take you outside $%
NP(n,m)$. So we must find sequences of profiles all of which are in $NP(n,m)$%
. In addition, we may have to vary the order of individuals whose
preferences are to be changed (an order which is fixed in standard sequence
arguments). To deal with these problems, we start with a definition and then
a lemma.\medskip

\qquad Given two profiles $u$ and $u^{\ast }$ in $NP(n,m)$, consider a
sequence of profiles $u_{1}$, $u_{2}$, ..., $u_{T}$, all in $NP(n,m)$, such
that\medskip

\qquad \qquad \qquad 1. $u_{1}=u$;\medskip

\qquad \qquad \qquad 2. $u_{T}=u^{\ast }$;\medskip

\qquad \qquad \qquad 3. For each $t$, $1\leq t<T$, there exists an $h\in N$

\qquad \qquad \qquad \qquad \qquad such that $u_{t}$ and $u_{t+1}$are $h$%
-variants;\medskip

\qquad \qquad \qquad 4. For each $t$, $1\leq t<T$, $u_{t}\mid S=u_{t+1}\mid
S $.\medskip

We call such a sequence an $S$\emph{-path} in $NP(n,m)$ from $u$ to $u^{\ast
}$. Note two obvious properties of $S$-paths:\medskip

\qquad 1. If $u_{1}$, $u_{2}$, ..., $u_{T}$ is an $S$-path from $u$ to $v$,
then $u_{T}$, $u_{T-1}$, ..., $u_{1}$ is an $S$-path from $v$ to $u$;\medskip

\qquad 2. If $u_{1}$, $u_{2}$, ..., $u_{T}$ is an $S$-path from $u$ to $v$,
and $v_{1}$, $v_{2}$, ..., $v_{T^{\ast }}$ is an $S$-path from $v$ to $w$,
then $u_{1}$, $u_{2}$, ..., $u_{T}$, $v_{1}$, $v_{2}$, ..., $v_{T^{\ast }}$
is an S-path from $u$ to $w$.\medskip

We next prove two lemmas regarding $S$-paths in $NP(n,m).$

\qquad \textbf{Lemma 3-1}. Suppose there exists an $x$ in $X$ such that
profiles $u$ and $u^{\ast }$ in $NP(n,m)$ agree on $S=X\backslash \{x\}$.
Then for each of $u$, $u^{\ast }$, there is as $S$-path in $NP(n,m)$ from $u$
to $u^{\ast }$

\textbf{Proof}:\medskip

We show that for each of $u$, $u^{\ast }$, there is as $S$-path in $NP(n,m)$
from $u$ to $u^{\ast \ast }$ \medskip 
\[
u^{\ast \ast }:%
\begin{tabular}{|c|c|c|c|}
\hline
\emph{1} & \emph{2} & $\cdots $ & \emph{n} \\ \hline
$x$ & $u(2)\mid S$ & $\cdots $ & $u(n)\mid S$ \\ 
$u(1)\mid S$ & $x$ &  & $x$ \\ \hline
\end{tabular}%
\]

Then an $S$-path from $u$ to $u^{\ast }$ in $NP(n,m)$ will be found by first
following the path from $u$ to $u^{\ast \ast }$ and then following - in
reverse - the path from $u^{\ast }$ to $u^{\ast \ast }$ (see properties 1
and 2 above).

\qquad If $x$ is individual \#1's top alternative, then create the path to $%
u^{\ast \ast }$ by taking $x$ down to the bottom for each $i>1$ in turn.
\medskip

\qquad Otherwise, we proceed by steps, raising $x$ one rank in $\#1$'s
ordering until it is $\#1$'s top. We describe one such step. Without loss of
generality, suppose $\#1$'s ordering, $1:a...bcx...$, with $x$ in the
position just below $c$. If anyone else has $c\succ _{u(j)}x$, we could
bring $x$ up just above $c$ in $\#1$'s ordering and still be in $NP(n,m)$.
So assume $x\succ _{u(j)}c$ for all $j>1$.\medskip 
\[
u:%
\begin{tabular}{|c|c|c|c|}
\hline
\emph{1} & \emph{2} & \emph{\ldots } & \emph{n} \\ \hline
$a$ & $\vdots $ &  & $\vdots $ \\ 
$\vdots $ & $x$ & $\ldots $ & $x$ \\ 
$b$ & $\vdots $ &  & $\vdots $ \\ 
$c$ & $c$ & $\ldots $ & $c$ \\ 
$x$ & $\vdots $ &  & $\vdots $ \\ 
$\vdots $ &  &  &  \\ \hline
\end{tabular}%
\]%
Whatever is between $x$ and $c$ (if anything) for $\#2$, is either below $x$
for $\#1$ or in $\{a,...,b\}$. But if we are in $NP(n,m)$, each of the
alternatives in $\{a,...,b\}$ above $c$ for $\#2$ must be below $c$, and so
below $x$, for some individual $j>2$. Then we can bring $x$ down just below $%
c$ for $\#2$ and stay in $NP(n,m)$. That allows us to bring $x$ up just
above $c$ in $\#1$'s ordering and still be in $NP(n,m)$. This continues
until x has been raised to 1's top. \ \ \ \ $\square $\medskip

Now we extend that result in the \emph{S-Path Lemma}. Note that this lemma
is \textbf{not} about a (strategy-proof) rule $g$ - it is solely about the
structure of $NP(n,m)$.\medskip

\textbf{Lemma 3-2 }For $S$ any subset of $X$, let $u$ and $u^{\ast }$ be two
profiles in $NP(n,m)$ with $u\mid S=u^{\ast }\mid S$. Then there exists an $%
S $-path in $NP(n,m)$ from $u$ to $u^{\ast }$.\medskip

\qquad Proof: To establish the existence of a sequence of profiles in $%
NP(n,m)$ from $u$ to $u^{\ast }$, we will show that there are $S$-path
sequences in $NP(n,m)$ from each of $u$ and $u^{\ast }$ to a profile $%
u^{\ast \ast }$ we describe shortly. Then an $S$-path from $u$ to $u^{\ast } 
$ in $NP(n,m)$ will be found by first following the path from $u$ to $%
u^{\ast \ast }$ and then following - in reverse - the path from $u^{\ast }$
to $u^{\ast \ast }$ (see properties 1 and 2 above). We suppose $X\backslash
S=\{x,y,...,z\}$ and let $xy...z$ be a fixed ordering on $X\backslash S$.
Then $u^{\ast \ast }$ is given by\medskip 
\[
u^{\ast \ast }:%
\begin{tabular}{|c|c|c|c|}
\hline
\emph{1} & \emph{2} & $\cdots $ & \emph{n} \\ \hline
$x$ &  &  &  \\ 
$y$ & $u(2)\mid S$ &  & $u(n)\mid S$ \\ 
$\vdots $ &  &  &  \\ 
$z$ &  &  &  \\ 
& $z$ &  & $z$ \\ 
$u(1)\mid S$ & $\vdots $ & $\cdots $ & $\vdots $ \\ 
& $y$ &  & $y$ \\ 
& $x$ &  & $x$ \\ \hline
\end{tabular}%
\]

So individual $\#1$ ranks everything in $X\backslash S$ above everything in $%
S$ and ranks the elements of $S$ the same way as they are ranked in $u(1)$.
Individuals $2,...,n$ rank everything in $S$ above everything in $%
X\backslash S$ and rank the elements of $S$ the same way they rank those
elements at $u$. Finally individuals $2,...,n$ rank the elements of $%
X\backslash S$ as $z...yx$, the opposite of their ordering by $\#1$. \medskip

\qquad We show that there is an $S$-path in $NP(n,m)$ from both $u$ and $%
u^{\ast }$ (which agree on $S$) to $u^{\ast \ast }$ by a series of
applications of Lemma 3-1. \ In this case with $X\backslash S=\{x,y,...,z\}$%
, we first set $S_{1}=X\backslash \{x\}$ and so $X\backslash S_{1}=\{x\}$.
Then, by Lemma 3-1, there is a path from $u$ to the following profile $u_{1}$%
:\medskip 
\[
u_{1}:%
\begin{tabular}{|c|c|c|c|}
\hline
\emph{1} & \emph{2} & $\cdots $ & \emph{n} \\ \hline
$x$ & $u(2)\mid S_{1}$ & $\cdots $ & $u(n)\mid S_{1}$ \\ 
$u(1)\mid S_{1}$ & $x$ &  & $x$ \\ \hline
\end{tabular}%
\]%
Next, take $X=S_{1}$ and $S_{2}=X\backslash \{y\}$. By Lemma 3-1, there is
an $S_{1}$-path and so an $S$-path from $u_{1}|S_{1}$ to the following
profile:\medskip

\[
u_{2}:%
\begin{tabular}{|c|c|c|c|}
\hline
\emph{1} & \emph{2} & $\cdots $ & \emph{n} \\ \hline
$y$ & $u_{1}(2)\mid S_{2}$ & $\cdots $ & $u_{1}(n)\mid S_{2}$ \\ 
$u_{1}(1)\mid S_{2}$ & $y$ &  & $y$ \\ \hline
\end{tabular}%
\medskip 
\]%
\medskip

If we take each profile in that path and insert an $x$ at the top of \#1's
ordering and $x$ at the bottom for everyone else, we get an $S$-path from $%
u_{1}$ to u$_{3_{\text{:}}}$

\[
u_{3}:%
\begin{tabular}{|c|c|c|c|}
\hline
\emph{1} & \emph{2} & \emph{\ldots } & \emph{n} \\ \hline
$x$ & $u_{1}(2)\mid S_{2}$ &  & $u_{1}(n)\mid S_{2}$ \\ 
$y$ & $y$ & $\ldots $ & $y$ \\ 
$u_{1}(1)\mid S_{2}$ & $x$ &  & $x$ \\ \hline
\end{tabular}%
\]%
Combining these two paths sequentially yields a path from u to u$_{3}$.
Continuing in this fashion we get an $S_{2}$-path and so an $S$-path from $u$
to $u^{\ast \ast }$. Continuing in this pattern yields a path from $u^{\ast
} $ to $u^{\ast \ast }$, which was our goal. \ \ \ \ \ $\square $\medskip

\medskip \qquad \textbf{Remark 1.} The $S$-path lemma does not require that $%
S$ be Range($g$) for a strategy-proof $g$, though that will be our primary
application, as in the Equivalence Theorem just below.\medskip

\qquad \textbf{Remark 2. }In the $S$-path lemma\textbf{,} $S$ can even be a
singleton, showing that there is a path in $NP(n,m)$ from any profile in $%
NP(n,m)$ to any other.\medskip

These two remarks will not be used in this paper, but they serve to
emphasize that the $S$-path lemma is \textbf{not} about a social choice
rule, but rather about the domain $NP(n,m)$.\medskip

\qquad \textbf{Equivalence Theorem 3-3.} For $m,n\geq 3$, and any rule $g$
that is strategy-proof on $NP(n,,m)$ and has range $S$, if $u\mid S=u^{\ast
}\mid S$, then%
\[
g(u)=g(u^{\ast })\text{.} 
\]%
\qquad \medskip \textbf{Proof.} Let $u_{1},u_{2},...,u_{T}$ be an $S$-path
in $NP(n,m)$ from $u$ to $u^{\ast }$ as guaranteed by the $S$-path Lemma. If 
$g(u)\neq g(u\ast )$, there must be a $t$, $1\leq t<T$, such that $%
g(u_{t})\neq g(u_{t+1})$, where $u_{t}$ and $u_{t+1}$ are $h$-variants. But
individual $h$ orders $g(u_{t})$ and $g(u_{t+1})$ the same at $u_{t}$ and $%
u_{t+1}$ since $u_{t}\mid S=u_{t+1}\mid S$, so $g$ must be manipulable by $h$
at either $u_{t}$ or $u_{t+1}$. Therefore, $g(u)=g(u^{\ast })$. \ \ \ \ $%
\square $\medskip

This theorem will be the key at the conclusion to extending our analysis
from full-range rules to the more general case.

\section{Induction basis: 3 individuals and 3 alternatives.\protect\medskip}

\qquad Now let $g$ be a strategy-proof social choice rule on $NP(n,m)$ with $%
\mid $Range($g$)$\mid $ $\geq 3$. We would like to show $g$ is dictatorial,
i.e., there exists an individual $i$ such that for each profile $u$ in $%
NP(n,m)$, alternative $g(u)$ is the top-ranked element in $u(i)$ restricted
to Range($g$). In this and the next two sections, we do this for the special
case where Range($g$) $=X$. Our goal, then, is to prove, for $m,n\geq 3$,
the proposition $SP(n,m)$\textbf{:} \textbf{every strategy-proof rule on }$%
NP(n,m)$\textbf{\ is dictatorial if it has full range.}\medskip

Our analysis proceeds by induction on $n$ first and then by induction on $m$%
. The basis step deals with $m=3$ and $n=3$.\medskip

\textbf{Theorem 4-3 (Basis). }$\mathbf{SP(3,3)}$\textbf{.} That is, let $g$
be a strategy-proof social choice rule on $NP(3,3)$. If Range($g$) $=X$,
then $g$ is dictatorial: there exists an individual $i$ such that for each
profile $u$ in $NP(3,3)$, $g(u)$ is the top-ranked element in $u(i)$.\medskip

\textbf{Proof:} Assume $X=\{a,b,c\}$ and $N=\{1,2,3\}$. \ The proof consists
of three steps.\medskip

\qquad \textbf{(Step 1)} For any strategy-proof rule $g$ on $NP(3,3)$, we
show how choice at a profile in $VP$, where $VP$ is the subdomain of $%
NP(3,3) $ consisting of voting paradox profiles (where each alternative
loses to some other under majority voting), will lead to individual
decisiveness for one alternative over another, on all of $NP(3,3)$.\medskip

\qquad \textbf{(Step 2) }Then we show that if the range of strategy-proof $g$%
, restricted to $VP$, is all of $X$, there must be a dictator for $g$%
.\medskip

\qquad \textbf{(Step 3)} Last, we show that if the range of $g$ on $NP(3,3)$
is $X$ then the range of $g$ restricted to $VP$ is also $X$.\medskip

\qquad We begin with Step 1, the decisiveness result.\medskip

\qquad \textbf{Lemma 4-4.} If rule $g$ is strategy-proof on $NP(3,3)$ then,
for every profile $u$ in $VP$, if alternative $\alpha \in X=\{a,b,c\}$ is
chosen at $u$, and $\alpha $ is top for individual $j$ while $\beta \in X$
is $j$'s bottom at $u$, then $g$ never chooses $\beta $ at any profile in $%
NP $ where $j$ prefers $\alpha $ to $\beta $; in this case we say $j$ is 
\textit{decisive} for $\alpha $ against $\beta $ on $NP$.\medskip

\qquad (If $j$ is decisive for $\alpha $ against $\beta $, then we write $%
\alpha D_{j}\beta $.)\medskip

\qquad \textbf{Proof:} Consider a voting paradox profile $u$:\medskip 
\[
u:%
\begin{tabular}{|c|c|c|}
\hline
\emph{1} & \emph{2} & \emph{3} \\ \hline
$a$ & $b$ & $c$ \\ 
$b$ & $c$ & $a$ \\ 
$c$ & $a$ & $b$ \\ \hline
\end{tabular}%
\]%
Without loss of generality, assume $g(u)=a$. We want to show for every
profile $v$ in $NP$ with $a\succ _{v(1)}c$ will have $g(v)\neq c$. Such
profiles have $v(1)=abc$, $acb$, or $bac$.\medskip

\qquad \textbf{Case 1}. $v(1)=abc$. So $v$ is:\medskip 
\[
v:%
\begin{tabular}{|c|c|c|}
\hline
\emph{1} & \emph{2} & \emph{3} \\ \hline
$a$ &  &  \\ 
$b$ & $v(2)$ & $v(3)$ \\ 
$c$ &  &  \\ \hline
\end{tabular}%
\]%
\medskip Note $g(u^{1})=a$ at $u^{1}$, a $2$-variant of $u$:\medskip 
\[
u^{1}:%
\begin{tabular}{|c|c|c|}
\hline
\emph{1} & \emph{2} & \emph{3} \\ \hline
$a$ & $c$ & $c$ \\ 
$b$ & $b$ & $a$ \\ 
$c$ & $a$ & $b$ \\ \hline
\end{tabular}%
\]%
\medskip or $\#2$ will manipulate to $u^{1}$ from $u$. Then look at $u^{2}$,
a $3$-variant of $u^{1}$:\medskip 
\[
u^{2}:%
\begin{tabular}{|c|c|c|}
\hline
\emph{1} & \emph{2} & \emph{3} \\ \hline
$a$ & $c$ &  \\ 
$b$ & $b$ & $v(3)$ \\ 
$c$ & $a$ &  \\ \hline
\end{tabular}%
\]%
\medskip If $g(u^{2})=c$, then $\#3$ would manipulate from $u^{1}$ to $u^{2}$%
. So $g(u^{2})\neq c$. But then $g(v)\neq c$ or $\#2$ would manipulate from $%
u^{2}$ to $v$. Therefore $c$ is not chosen at any profile in $NP$ at which
individual $\#1$'s ordering is $abc$.\medskip

\qquad \textbf{Case 2}. $v(1)=acb$. One profile with person $\#1$ having
ordering $acb$ is $u^{^{\prime }}$, a simple $1$-variant of $u$
above:\medskip 
\[
u^{\prime }:%
\begin{tabular}{|c|c|c|}
\hline
\emph{1} & \emph{2} & \emph{3} \\ \hline
$a$ & $b$ & $c$ \\ 
$c$ & $c$ & $a$ \\ 
$b$ & $a$ & $b$ \\ \hline
\end{tabular}%
\]%
\medskip $g(u^{\prime })=a$ or $\#1$ would manipulate to $u$ with $u(1)=abc$%
.\medskip\ 

\qquad We trace out the consequences of this by considering the next figure,
where individual $\#1$ has ordering $acb$, $\#2$'s possible orderings are in
different rows, and $\#3$'s preference orderings are in different columns.
[Black cells are not in $NP(3,3)$; gray background cells are in $VP$.
\textquotedblleft $\lnot \ c$\textquotedblright\ in a cell means $c$ is not
chosen there.] There is an \textquotedblleft $a$\textquotedblright\ in cell $%
(4,5)$, i.e., row $4$ and column $5$, in the main body of the table (i.e.,
not including the preference label row and preference label column)
indicating the outcome at that cell, which is the profile $u^{\prime }$.
There is an \textquotedblleft $a$\textquotedblright\ above that in $(3,5)$
or $\#2$ (the row player) would manipulate from $(4,5)$ to $(3,5)$.
Similarly, there are \textquotedblleft $a$\textquotedblright s in cells $%
(4,1)$ and $(4,2)$ or $\#3$ (the column player) would manipulate from there
to $(4,5)$. An \textquotedblleft $a$\textquotedblright\ appears in $(6,1)$
or $\#2$ would manipulate there from $(4,1)$. Three cells in row $4$ are
labeled \textquotedblleft $\lnot \ c$\textquotedblright\ (i.e.,
\textquotedblleft not $c$\textquotedblright ) because $c$ in such a cell
would lead $\#3$ to manipulate there from $(4,5)$. Similarly for two cells
in row $3$.\medskip

\textbf{1: acb}\medskip

\begin{tabular}{|l|c|c|c|c|c|c|}
\hline
$%
\begin{array}{c}
\text{ \ \ \ \ \ }3 \\ 
\\ 
2\text{ \ \ \ \ \ }%
\end{array}%
$ & $%
\begin{tabular}{l}
$a$ \\ 
$b$ \\ 
$c$%
\end{tabular}%
$ & $%
\begin{tabular}{l}
$a$ \\ 
$c$ \\ 
$b$%
\end{tabular}%
$ & $%
\begin{tabular}{l}
$b$ \\ 
$a$ \\ 
$c$%
\end{tabular}%
$ & $%
\begin{tabular}{l}
$b$ \\ 
$c$ \\ 
$a$%
\end{tabular}%
$ & $%
\begin{tabular}{l}
$c$ \\ 
$a$ \\ 
$b$%
\end{tabular}%
$ & $%
\begin{tabular}{l}
$c$ \\ 
$b$ \\ 
$a$%
\end{tabular}%
$ \\ \hline
$%
\begin{tabular}{l}
\\ 
$a$ $b$ $c$ \\ 
\begin{tabular}{l}
\end{tabular}%
\end{tabular}%
$ & $\ \ \ \ 
\cellcolor[gray]{0.2}%
\ \ \ \ $ & $\ \ \ \ 
\cellcolor[gray]{0.2}%
\ \ \ \ \ $ & $\ \ \ \ \ 
\cellcolor[gray]{0.2}%
\ \ \ \ $ & $\ \ \ \ \ \ \ \ \ $ & $\ \ \ \ 
\cellcolor[gray]{0.2}%
\ \ \ \ \ $ & $\ \ \ \ \ \ \ \ \ $ \\ \hline
$%
\begin{tabular}{l}
\\ 
$a$ $c$ $b$ \\ 
\begin{tabular}{l}
\end{tabular}%
\end{tabular}%
$ & 
\cellcolor[gray]{0.2}
& 
\cellcolor[gray]{0.2}
& 
\cellcolor[gray]{0.2}
&  & 
\cellcolor[gray]{0.2}
& 
\cellcolor[gray]{0.2}
\\ \hline
$%
\begin{tabular}{l}
\\ 
$b$ $a\ c$ \\ 
\begin{tabular}{l}
\end{tabular}%
\end{tabular}%
$ & 
\cellcolor[gray]{0.2}
& 
\cellcolor[gray]{0.2}
& 
\cellcolor[gray]{0.2}
& $%
\begin{tabular}{l}
\\ 
\\ 
$\lnot c$%
\end{tabular}%
$ & $a$ & 
\cellcolor[gray]{0.75}%
$%
\begin{tabular}{l}
\\ 
\\ 
$\lnot c$%
\end{tabular}%
$ \\ \hline
$%
\begin{tabular}{l}
\\ 
$b\ c\ a$ \\ 
\begin{tabular}{l}
\end{tabular}%
\end{tabular}%
$ & $a$ & $a$ & $%
\begin{tabular}{l}
\\ 
\\ 
$\lnot c$%
\end{tabular}%
$ & $%
\begin{tabular}{l}
\\ 
\\ 
$\lnot c$%
\end{tabular}%
$ & $a$ & $%
\begin{tabular}{l}
\\ 
\\ 
$\lnot c$%
\end{tabular}%
$ \\ \hline
$%
\begin{tabular}{l}
\\ 
$c\ a$ $b$ \\ 
\begin{tabular}{l}
\end{tabular}%
\end{tabular}%
$ & 
\cellcolor[gray]{0.2}
& 
\cellcolor[gray]{0.2}
&  &  & 
\cellcolor[gray]{0.2}
& 
\cellcolor[gray]{0.2}
\\ \hline
$%
\begin{tabular}{l}
\\ 
$c$ $b$ $a$ \\ 
\begin{tabular}{l}
\end{tabular}%
\end{tabular}%
$ & $a$ & 
\cellcolor[gray]{0.2}
& 
\cellcolor[gray]{0.75}
&  & 
\cellcolor[gray]{0.2}
& 
\cellcolor[gray]{0.2}
\\ \hline
\end{tabular}%
\medskip \medskip

\qquad We extend this analysis in the next table, where the results just
obtained are emboldened and in regular (non-Italic). Now $c$ can't be chosen
at $(6,3)$ or $\#3$ will manipulate to $(6,1)$. And that means $c$ isn't
chosen at $(5,3)$ or $\#2$ would manipulate from $(6,3)$ to $(5,3)$. At
cells $(1,4)$ and $(1,6)$, alternative $c$ won't be chosen because $c$ is $%
\#2$'s worst there and $\#2$ would manipulate to row $3$ (or $4$).\medskip
\medskip \pagebreak

\textbf{1: acb}\medskip

\begin{tabular}{|l|c|c|c|c|c|c|}
\hline
$%
\begin{array}{c}
\text{ \ \ \ }3 \\ 
\\ 
2\text{ \ \ \ \ \ }%
\end{array}%
$ & $%
\begin{tabular}{l}
$a$ \\ 
$b$ \\ 
$c$%
\end{tabular}%
$ & $%
\begin{tabular}{l}
$a$ \\ 
$c$ \\ 
$b$%
\end{tabular}%
$ & $%
\begin{tabular}{l}
$b$ \\ 
$a$ \\ 
$c$%
\end{tabular}%
$ & $%
\begin{tabular}{l}
$b$ \\ 
$c$ \\ 
$a$%
\end{tabular}%
$ & $%
\begin{tabular}{l}
$c$ \\ 
$a$ \\ 
$b$%
\end{tabular}%
$ & $%
\begin{tabular}{l}
$c$ \\ 
$b$ \\ 
$a$%
\end{tabular}%
$ \\ \hline
$%
\begin{tabular}{l}
\\ 
$a$ $b$ $c$ \\ 
\begin{tabular}{l}
\end{tabular}%
\end{tabular}%
$ & $\ \ \ \ 
\cellcolor[gray]{0.2}%
\ \ \ \ $ & $\ \ \ \ 
\cellcolor[gray]{0.2}%
\ \ \ \ \ $ & $\ \ \ \ \ 
\cellcolor[gray]{0.2}%
\ \ \ \ $ & $\ \ \ \ 
\begin{tabular}{l}
\\ 
\\ 
$\lnot c$%
\end{tabular}%
\ \ \ \ \ $ & $\ \ \ \ 
\cellcolor[gray]{0.2}%
\ \ \ \ \ $ & $\ \ \ \ 
\begin{tabular}{l}
\\ 
\\ 
$\lnot c$%
\end{tabular}%
\ \ \ \ \ $ \\ \hline
$%
\begin{tabular}{l}
\\ 
$a$ $c$ $b$ \\ 
\begin{tabular}{l}
\end{tabular}%
\end{tabular}%
$ & 
\cellcolor[gray]{0.2}
& 
\cellcolor[gray]{0.2}
& 
\cellcolor[gray]{0.2}
&  & 
\cellcolor[gray]{0.2}
& 
\cellcolor[gray]{0.2}
\\ \hline
$%
\begin{tabular}{l}
\\ 
$b$ $a\ c$ \\ 
\begin{tabular}{l}
\end{tabular}%
\end{tabular}%
$ & 
\cellcolor[gray]{0.2}
& 
\cellcolor[gray]{0.2}
& 
\cellcolor[gray]{0.2}
& $%
\begin{tabular}{l}
\\ 
\\ 
$\lnot \mathbf{c}$%
\end{tabular}%
$ & $\mathbf{a}$ & 
\cellcolor[gray]{0.75}%
$%
\begin{tabular}{l}
\\ 
\\ 
$\lnot \mathbf{c}$%
\end{tabular}%
$ \\ \hline
$%
\begin{tabular}{l}
\\ 
$b\ c\ a$ \\ 
\begin{tabular}{l}
\end{tabular}%
\end{tabular}%
$ & $\mathbf{a}$ & $\mathbf{a}$ & $%
\begin{tabular}{l}
\\ 
\\ 
$\lnot \mathbf{c}$%
\end{tabular}%
$ & $%
\begin{tabular}{l}
\\ 
\\ 
$\lnot \mathbf{c}$%
\end{tabular}%
$ & $\mathbf{a}$ & $%
\begin{tabular}{l}
\\ 
\\ 
$\lnot \mathbf{c}$%
\end{tabular}%
$ \\ \hline
$%
\begin{tabular}{l}
\\ 
$c\ a$ $b$ \\ 
\begin{tabular}{l}
\end{tabular}%
\end{tabular}%
$ & 
\cellcolor[gray]{0.2}
& 
\cellcolor[gray]{0.2}
& $%
\begin{tabular}{l}
\\ 
\\ 
$\lnot c$%
\end{tabular}%
$ &  & 
\cellcolor[gray]{0.2}
& 
\cellcolor[gray]{0.2}
\\ \hline
$%
\begin{tabular}{l}
\\ 
$c$ $b$ $a$ \\ 
\begin{tabular}{l}
\end{tabular}%
\end{tabular}%
$ & $\mathbf{a}$ & 
\cellcolor[gray]{0.2}
& 
\cellcolor[gray]{0.75}%
$%
\begin{tabular}{l}
\\ 
\\ 
$\lnot c$%
\end{tabular}%
$ &  & 
\cellcolor[gray]{0.2}
& 
\cellcolor[gray]{0.2}
\\ \hline
\end{tabular}%
\medskip \bigskip

\qquad In the next table, again all the preceding results are emboldened .
Suppose $c$ is chosen at cell $(6,4)$ indicated by \textquotedblleft $c?$%
\textquotedblright . Then $b$ is not chosen in $(6,3)$ or $\#3$ will
manipulate there from $(6,4)$. Therefore $a$ is chosen in $(6,3)$. In turn,
since $a$ is $\#2$'s bottom at $(6,3)$, $a$ must also be chosen at $(5,3)$
and $(4,3)$. Since $c$ is chosen at $(6,4)$, $a$ can't be chosen at $(4,4)$,
or $\#2$ would manipulate from there to $(6,4)$. Hence $b$ is chosen at $%
(4,4)$. But that leads to a contradiction, since now $\#3$ would manipulate
from $(4,3)$ to $(4,4)$. From all that, we can conclude $c$ is not chosen at 
$(6,4)$ and that implies $c$ is chosen nowhere in column $4$ or $\#2$ would
manipulate up from $(6,4)$. We have determined that $c$ is chosen in none of
the $NP(3,3)$ profiles with individual $\#1$'s ordering $acb$.\medskip
\medskip \pagebreak

\textbf{1: acb}\medskip

\begin{tabular}{|l|c|c|c|c|c|c|}
\hline
$%
\begin{array}{c}
\text{ \ \ \ }3 \\ 
\\ 
2\text{ \ \ \ \ \ }%
\end{array}%
$ & $%
\begin{tabular}{l}
$a$ \\ 
$b$ \\ 
$c$%
\end{tabular}%
$ & $%
\begin{tabular}{l}
$a$ \\ 
$c$ \\ 
$b$%
\end{tabular}%
$ & $%
\begin{tabular}{l}
$b$ \\ 
$a$ \\ 
$c$%
\end{tabular}%
$ & $%
\begin{tabular}{l}
$b$ \\ 
$c$ \\ 
$a$%
\end{tabular}%
$ & $%
\begin{tabular}{l}
$c$ \\ 
$a$ \\ 
$b$%
\end{tabular}%
$ &  \\ \hline
$%
\begin{tabular}{l}
\\ 
$a$ $b$ $c$ \\ 
\begin{tabular}{l}
\end{tabular}%
\end{tabular}%
$ & $\ \ \ \ 
\cellcolor[gray]{0.2}%
\ \ \ \ $ & $\ \ \ \ 
\cellcolor[gray]{0.2}%
\ \ \ \ \ $ & $\ \ \ \ \ 
\cellcolor[gray]{0.2}%
\ \ \ \ $ & $\ \ \ \ 
\begin{tabular}{l}
\\ 
\\ 
$\lnot \mathbf{c}$%
\end{tabular}%
\ \ \ \ \ $ & $\ \ \ \ 
\cellcolor[gray]{0.2}%
\ \ \ \ \ $ & $\ \ \ \ 
\begin{tabular}{l}
\\ 
\\ 
$\lnot \mathbf{c}$%
\end{tabular}%
\ \ \ \ \ $ \\ \hline
$%
\begin{tabular}{l}
\\ 
$a$ $c$ $b$ \\ 
\begin{tabular}{l}
\end{tabular}%
\end{tabular}%
$ & 
\cellcolor[gray]{0.2}
& 
\cellcolor[gray]{0.2}
& 
\cellcolor[gray]{0.2}
&  & 
\cellcolor[gray]{0.2}
& 
\cellcolor[gray]{0.2}
\\ \hline
$%
\begin{tabular}{l}
\\ 
$b$ $a\ c$ \\ 
\begin{tabular}{l}
\end{tabular}%
\end{tabular}%
$ & 
\cellcolor[gray]{0.2}
& 
\cellcolor[gray]{0.2}
& 
\cellcolor[gray]{0.2}
& $%
\begin{tabular}{l}
\\ 
\\ 
$\lnot \mathbf{c}$%
\end{tabular}%
$ & $\mathbf{a}$ & 
\cellcolor[gray]{0.75}%
$%
\begin{tabular}{l}
\\ 
\\ 
$\lnot \mathbf{c}$%
\end{tabular}%
$ \\ \hline
$%
\begin{tabular}{l}
\\ 
$b\ c\ a$ \\ 
\begin{tabular}{l}
\end{tabular}%
\end{tabular}%
$ & $\mathbf{a}$ & $\mathbf{a}$ & $%
\begin{tabular}{l}
\\ 
$a$ \\ 
$\lnot \mathbf{c}$%
\end{tabular}%
$ & $%
\begin{tabular}{l}
\\ 
$b$ \\ 
$\lnot \mathbf{c,}\lnot a$%
\end{tabular}%
$ & $\mathbf{a}$ & $%
\begin{tabular}{l}
\\ 
\\ 
$\lnot \mathbf{c}$%
\end{tabular}%
$ \\ \hline
$%
\begin{tabular}{l}
\\ 
$c\ a$ $b$ \\ 
\begin{tabular}{l}
\end{tabular}%
\end{tabular}%
$ & 
\cellcolor[gray]{0.2}
& 
\cellcolor[gray]{0.2}
& $%
\begin{tabular}{l}
\\ 
$a$ \\ 
$\lnot \mathbf{c}$%
\end{tabular}%
$ &  & 
\cellcolor[gray]{0.2}
& 
\cellcolor[gray]{0.2}
\\ \hline
$%
\begin{tabular}{l}
\\ 
$c$ $b$ $a$ \\ 
\begin{tabular}{l}
\end{tabular}%
\end{tabular}%
$ & $\mathbf{a}$ & 
\cellcolor[gray]{0.2}
& 
\cellcolor[gray]{0.75}%
$%
\begin{tabular}{l}
\\ 
$a$ \\ 
$\lnot \mathbf{c,}\lnot b$%
\end{tabular}%
$ & $c?$ & 
\cellcolor[gray]{0.2}
& 
\cellcolor[gray]{0.2}
\\ \hline
\end{tabular}%
\medskip \bigskip \medskip

\qquad \textbf{Case 3}. $v(1)=bac$. In the next figure, we show the $NP(3,3)$
profiles for which individual $\#1$ has ordering $bac$. In twelve cells, we
have entered \textquotedblleft $\lnot \ c$\textquotedblright\ because, were $%
c$ to be chosen at one of those cells, $\#1$ (for whom $c$ is worst) would
manipulate by changing his ordering to $abc$ or $acb$, whichever would yield
a profile still in $NP(3,3)$ because we know that $c$ wouldn't be chosen at
the resulting profile.\medskip \pagebreak \medskip

\textbf{1: bac}\medskip

\begin{tabular}{|l|c|c|c|c|c|c|}
\hline
$%
\begin{array}{c}
\text{ \ \ \ }3 \\ 
\\ 
2\text{ \ \ \ \ \ }%
\end{array}%
$ & $%
\begin{tabular}{l}
$a$ \\ 
$b$ \\ 
$c$%
\end{tabular}%
$ & $%
\begin{tabular}{l}
$a$ \\ 
$c$ \\ 
$b$%
\end{tabular}%
$ & $%
\begin{tabular}{l}
$b$ \\ 
$a$ \\ 
$c$%
\end{tabular}%
$ & $%
\begin{tabular}{l}
$b$ \\ 
$c$ \\ 
$a$%
\end{tabular}%
$ & $%
\begin{tabular}{l}
$c$ \\ 
$a$ \\ 
$b$%
\end{tabular}%
$ & $%
\begin{tabular}{l}
$c$ \\ 
$b$ \\ 
$a$%
\end{tabular}%
$ \\ \hline
$%
\begin{tabular}{l}
\\ 
$a$ $b$ $c$ \\ 
\begin{tabular}{l}
\end{tabular}%
\end{tabular}%
$ & $\ \ \ \ 
\cellcolor[gray]{0.2}%
\ \ \ \ $ & $\ \ \ \ 
\cellcolor[gray]{0.2}%
\ \ \ \ \ $ & $\ \ \ \ \ 
\cellcolor[gray]{0.2}%
\ \ \ \ $ & $\ \ \ \ 
\cellcolor[gray]{0.2}%
\ \ \ \ \ $ & $\ \ \ \ \ \ \ \ \ $ & $%
\begin{tabular}{l}
\\ 
\\ 
$\lnot c$%
\end{tabular}%
$ \\ \hline
$%
\begin{tabular}{l}
\\ 
$a$ $c$ $b$ \\ 
\begin{tabular}{l}
\end{tabular}%
\end{tabular}%
$ & 
\cellcolor[gray]{0.2}
& 
\cellcolor[gray]{0.2}
& 
\cellcolor[gray]{0.2}
& $%
\begin{tabular}{l}
\\ 
\\ 
$\lnot c$%
\end{tabular}%
$ &  & 
\cellcolor[gray]{0.80}%
$%
\begin{tabular}{l}
\\ 
\\ 
$\lnot c$%
\end{tabular}%
$ \\ \hline
$%
\begin{tabular}{l}
\\ 
$b$ $a\ c$ \\ 
\begin{tabular}{l}
\end{tabular}%
\end{tabular}%
$ & 
\cellcolor[gray]{0.2}
& 
\cellcolor[gray]{0.2}
& 
\cellcolor[gray]{0.2}
& 
\cellcolor[gray]{0.2}
& $%
\begin{tabular}{l}
\\ 
\\ 
$\lnot c$%
\end{tabular}%
$ & 
\cellcolor[gray]{0.75}
\\ \hline
$%
\begin{tabular}{l}
\\ 
$b\ c\ a$ \\ 
\begin{tabular}{l}
\end{tabular}%
\end{tabular}%
$ & 
\cellcolor[gray]{0.2}
& $%
\begin{tabular}{l}
\\ 
\\ 
$\lnot c$%
\end{tabular}%
$ & 
\cellcolor[gray]{0.2}
& 
\cellcolor[gray]{0.2}
& $%
\begin{tabular}{l}
\\ 
\\ 
$\lnot c$%
\end{tabular}%
$ & 
\cellcolor[gray]{0.2}
\\ \hline
$%
\begin{tabular}{l}
\\ 
$c\ a$ $b$ \\ 
\begin{tabular}{l}
\end{tabular}%
\end{tabular}%
$ & 
\cellcolor[gray]{0.2}
& 
\cellcolor[gray]{0.2}
& $%
\begin{tabular}{l}
\\ 
\\ 
$\lnot c$%
\end{tabular}%
$ & $%
\begin{tabular}{l}
\\ 
\\ 
$\lnot c$%
\end{tabular}%
$ &  & $%
\begin{tabular}{l}
\\ 
\\ 
$\lnot c$%
\end{tabular}%
$ \\ \hline
$%
\begin{tabular}{l}
\\ 
$c$ $b$ $a$ \\ 
\begin{tabular}{l}
\end{tabular}%
\end{tabular}%
$ & $%
\begin{tabular}{l}
\\ 
\\ 
$\lnot c$%
\end{tabular}%
$ & 
\cellcolor[gray]{0.80}%
$%
\begin{tabular}{l}
\\ 
\\ 
$\lnot c$%
\end{tabular}%
$ & 
\cellcolor[gray]{0.75}
& 
\cellcolor[gray]{0.2}
& $%
\begin{tabular}{l}
\\ 
\\ 
$\lnot c$%
\end{tabular}%
$ & 
\cellcolor[gray]{0.2}
\\ \hline
\end{tabular}%
\bigskip

\qquad Since $c$ is not chosen in cell $(5,6)$ where $c$ is at the top of $%
\#3$'s ordering, $c$ must not be chosen anywhere in row $5$ or else $\#3$
would manipulate there from $(5,6)$. But in row $5$, $c$ is at the top of $%
\#2$'s ordering and so $c$ can't be chosen anywhere in the table or $\#2$
would manipulate there from row $5$.\qquad $\square $\medskip

\qquad Step 2 of the proof consists of proving the following
corollaries.\medskip

\qquad \textbf{Corollary 4-4-1}. Suppose $g$ is a strategy-proof rule on $%
NP(3,3)$: then if $g$ is dictatorial on $VP$, it is dictatorial on all of $%
NP(3,3)$.\medskip

\qquad \textbf{Corollary 4-4-2}. Suppose $g$ is a strategy-proof rule on $%
NP(3,3)$: then if $g$ has singleton range on $VP$, it has singleton range on
all of $NP(3,3)$.\medskip

Proof of these first two corollaries is easy.\medskip

\qquad \textbf{Corollary 4-4-3}. Suppose $g$ is a strategy-proof rule on $%
NP(3,3)$: then:if $g$ has range $X$ on $VP,$ then it is dictatorial on all
of $NP(3,3).$\medskip

\qquad \textbf{Proof of Corollary 4-4-3}: It suffices, by Corollary 4-4-1,
to show that if $g$ has range $X$ on $VP$, then it is dictatorial on $VP$%
.\medskip

So consider the twelve profiles in $VP$:\medskip 
\[
\begin{tabular}{lllll}
\underline{$\mathbf{1}$} &  & \underline{$\mathbf{5}$} &  & \underline{$%
\mathbf{9}$} \\ 
&  &  &  &  \\ 
$\ 1:abc$ &  & $\ 1:bac$ &  & $\ 1:cab$ \\ 
$\ 2:bca\ \ \ \ \ a$ &  & $\ 2:acb\ \ \ \ \ not$ $c$ &  & $\ 2:abc$ \\ 
$\ 3:cab$ &  & $\ 3:cba$ &  & $\ 3:bca$ \\ 
&  &  &  &  \\ 
\underline{$\mathbf{2}$} &  & \underline{$\mathbf{6}$} &  & \underline{$%
\mathbf{10}$} \\ 
&  &  &  &  \\ 
$\ 1:abc$ &  & $\ 1:bac$ &  & $\ 1:cab$ \\ 
$\ 2:cab\ \ \ \ \ not$ $c$ &  & $\ 2:cba\ \ \ \ \ not$ $c$ &  & $\ 2:bca\ \
\ \ \ not$ $b$ \\ 
$\ 3:bca$ &  & $\ 3:acb$ &  & $\ 3:abc$ \\ 
&  &  &  &  \\ 
\underline{$\mathbf{3}$} &  & \underline{$\mathbf{7}$} &  & \underline{$%
\mathbf{11}$} \\ 
&  &  &  &  \\ 
$\ 1:acb$ &  & $\ 1:bca$ &  & $\ 1:cba$ \\ 
$\ 2:cba\ \ \ \ \ not$ $c$ &  & $\ 2:cab$ &  & $\ 2:bac$ \\ 
$\ 3:bac$ &  & $\ 3:abc$ &  & $\ 3:acb$ \\ 
&  &  &  &  \\ 
\underline{$\mathbf{4}$} &  & \underline{$\mathbf{8}$} &  & \underline{$%
\mathbf{12}$} \\ 
&  &  &  &  \\ 
$\ 1:acb$ &  & $\ 1:bca$ &  & $\ 1:cba$ \\ 
$\ 2:bac\ \ \ \ \ not$ $c$ &  & $\ 2:abc$ &  & $\ 2:acb$ \\ 
$\ 3:cba$ &  & $\ 3:cab$ &  & $\ 3:bac$%
\end{tabular}%
\]%
Without loss of generality, we assume that $g$ yields $a$ at the first
profile. Therefore, by Lemma 4, individual $\#1$ is decisive for $a$ against 
$c$, which we write $aD_{1}c$. That implies $c$ is not chosen at profiles 2,
3, 4, 5, and 6. Also $g$ can not yield $b$ at profile 10 or else $bD_{2}a$,
contrary to the outcome at profile 1.\pagebreak

\bigskip

\qquad \textbf{Profile 2}. Suppose b were chosen at profile 2:\medskip 
\[
\begin{tabular}{lllll}
$\mathbf{1}$ &  & $\mathbf{5}$ &  & $\mathbf{9}$ \\ 
&  &  &  &  \\ 
$\ 1:abc$ &  & $\ 1:bac$ &  & $\ 1:cab$ \\ 
$\ 2:bca\ \ \ \ \ \mathbf{a}$ &  & $\ 2:acb\ \ \ \ \ $\textbf{not} $\mathbf{c%
}$ &  & $\ 2:abc$ \ \ \ \ \textit{not} $a$ \\ 
$\ 3:cab$ &  & $\ 3:cba$ \ \ \ \ \textit{not} $a;$ \ $\therefore b$ &  & $\
3:bca$ \\ 
&  &  &  &  \\ 
$\mathbf{2}$ &  & $\mathbf{6}$ &  & $\mathbf{10}$ \\ 
&  &  &  &  \\ 
$\ 1:abc$ &  & $\ 1:bac$ &  & $\ 1:cab$ \\ 
$\ 2:cab\ \ \ \ \ $\textbf{not} $\mathbf{c};\ b?$ &  & $\ 2:cba\ \ \ \ \ $%
\textbf{not} $\mathbf{c}$ &  & $\ 2:bca\ \ \ \ \ $\textbf{not} $\mathbf{b}$
\\ 
$\ 3:bca$ &  & $\ 3:acb$ &  & $\ 3:abc$ \\ 
&  &  &  &  \\ 
$\mathbf{3}$ &  & $\mathbf{7}$ &  & $\mathbf{11}$ \\ 
&  &  &  &  \\ 
$\ 1:acb$ &  & $\ 1:bca$ &  & $\ 1:cba$ \\ 
$\ 2:cba\ \ \ \ \ $\textbf{not} $\mathbf{c};$ &  & $\ 2:cab$ &  & $\ 2:bac$
\\ 
$\ 3:bac\ \ \ \ \ not$ $a;\ \ \therefore b$ &  & $\ 3:abc$ &  & $\ 3:acb$ \\ 
&  &  &  &  \\ 
$\mathbf{4}$ &  & $\mathbf{8}$ &  & $\mathbf{12}$ \\ 
&  &  &  &  \\ 
$\ 1:acb$ &  & $\ 1:bca$ &  & $\ 1:cba$ \\ 
$\ 2:bac\ \ \ \ \ $\textbf{not} $\mathbf{c}$ &  & $\ 2:abc$ &  & $\ 2:acb$ \
\ \ \ \textit{not} $a$ \\ 
$\ 3:cba$ \ \ \ \ $not$ $a;$ $\ \therefore b$ &  & $\ 3:cab$ &  & $\ 3:bac$%
\end{tabular}%
\]%
Then by Lemma 4, $bD_{3}a$, which implies $a$ is not chosen at profiles 3,
4, 5, 9, and 12. But that implies $b$ is chosen at each of profiles 3, 4,
and 5. These results imply $bD_{3}c$, $bD_{2}c$, and $bD_{1}c$,
respectively. But at any profile in NP, at least one of the individuals must
prefer $b$ to $c$, so $c$ would never be chosen in $NP(3,3)$. So now we may
assume $b$ is not chosen at profile 2. Since we already know $c$ isn't
chosen there, $a$ must be chosen. Now $b$ can't be chosen at profile 9 or $%
bD_{3}a$, contrary to the choice of $a$ at profile 2. All results up to this
point are emboldened in the next display.\medskip

\pagebreak

\bigskip

\qquad \textbf{Profiles 3 and 4.} Suppose b were chosen at profile 3. Then
bD3c, and so c is not chosen at profiles 7, 9, 10, and 12.\medskip 
\[
\begin{tabular}{lllll}
$\mathbf{1}$ &  & $\mathbf{5}$ &  & $\mathbf{9}$ \\ 
&  &  &  &  \\ 
$\ 1:abc$ &  & $\ 1:bac$ &  & $\ 1:cab$ \\ 
$\ 2:bca\ \ \ \ \ \mathbf{a}$ &  & $\ 2:acb\ \ \ \ \ $\textbf{not} $\mathbf{c%
}$ &  & $\ 2:abc$ \ \ \ \ \textbf{not} $\mathbf{b}$ \\ 
$\ 3:cab$ &  & $\ 3:cba$ &  & $\ 3:bca$ \ \ \ \ $not$ $c;$ $\ \therefore a$
\\ 
&  &  &  &  \\ 
$\mathbf{2}$ &  & $\mathbf{6}$ &  & $\mathbf{10}$ \\ 
&  &  &  &  \\ 
$\ 1:abc$ &  & $\ 1:bac$ &  & $\ 1:cab$ \\ 
$\ 2:cab\ \ \ \ \ \mathbf{a}$ &  & $\ 2:cba\ \ \ \ \ $\textbf{not} $\mathbf{c%
}$ &  & $\ 2:bca\ \ \ \ \ $\textbf{not} $\mathbf{b};$ \\ 
$\ 3:bca$ &  & $\ 3:acb$ &  & $\ 3:abc$ \ \ \ \ $not$ $c;$ $\ \therefore a$
\\ 
&  &  &  &  \\ 
$\mathbf{3}$ &  & $\mathbf{7}$ &  & $\mathbf{11}$ \\ 
&  &  &  &  \\ 
$\ 1:acb$ &  & $\ 1:bca$ &  & $\ 1:cba$ \\ 
$\ 2:cba\ \ \ \ \ $\textbf{not} $\mathbf{c};$ \ \ $b?$ &  & $\ 2:cab$ \ \ \
\ $not$ $c$ &  & $\ 2:bac$ \ \ \ \ $not$ $c$ \\ 
$\ 3:bac$ &  & $\ 3:abc$ &  & $\ 3:acb$ \\ 
&  &  &  &  \\ 
$\mathbf{4}$ &  & $\mathbf{8}$ &  & $\mathbf{12}$ \\ 
&  &  &  &  \\ 
$\ 1:acb$ &  & $\ 1:bca$ &  & $\ 1:cba$ \\ 
$\ 2:bac\ \ \ \ \ $\textbf{not} $\mathbf{c}$ &  & $\ 2:abc$ \ \ \ \ $not$ $c$
&  & $\ 2:acb$ \ \ \ \ $not$ $c$ \\ 
$\ 3:cba$ &  & $\ 3:cab$ &  & $\ 3:bac$%
\end{tabular}%
\]%
At profile 9, neither $b$ nor $c$ is chosen, so a must be. Therefore, $%
aD_{2}c$ which implies $c$ is not chosen at profiles 8 or 11. Then $c$ would
not be in the range of $g$ restricted to $VP$. All that was a consequence of
assuming $b$ was chosen at profile 3. So $b$ isn't chosen there which means $%
a$ is. But then $aD_{1}b$, so $b$ is not chosen at profile 4, and therefore $%
a$ must be chosen there.\medskip \pagebreak

\bigskip

\textbf{Profile 5.\medskip }

\qquad Assume that $a$ is also chosen at profile 5. Then $aD_{2}b$ and so $b$
is not chosen at profiles 7, 8, and 12.\medskip

\qquad If $a$ were also chosen at profile 6, everyone would be decisive for $%
a$ against $b$ and $b$ would never be chosen on $NP(3,3)$. So $a$ is not
chosen there which means that $b$ is chosen there. But then $bD_{1}c$ and so 
$c$ is not chosen at profiles 7 and 8. This means $a$ is chosen at each of
those profiles and that implies everyone is decisive for $a$ against $c$ and 
$c$ would never be chosen on $NP(3,3)$. Therefore our assumption that $a$ is
chosen at profile 5 fails and $b$ must be chosen there. Since $b$ is chosen
at profile 5, $bD_{1}c$ and thus $c$ is not chosen at profiles 7 and
8.\medskip 
\[
\begin{tabular}{lllll}
$\mathbf{1}$ &  & $\mathbf{5}$ &  & $\mathbf{9}$ \\ 
&  &  &  &  \\ 
$\ 1:abc$ &  & $\ 1:bac$ &  & $\ 1:cab$ \\ 
$\ 2:bca\ \ \ \ \ \mathbf{a}$ &  & $\ 2:acb\ \ \ \ \ $\textbf{not} $\mathbf{c%
};$ \ $a?$ &  & $\ 2:abc$ \ \ \ \ \textbf{not} $\mathbf{b}$ \\ 
$\ 3:cab$ &  & $\ 3:cba$ &  & $\ 3:bca$ \\ 
&  &  &  &  \\ 
$\mathbf{2}$ &  & $\mathbf{6}$ &  & $\mathbf{10}$ \\ 
&  &  &  &  \\ 
$\ 1:abc$ &  & $\ 1:bac$ &  & $\ 1:cab$ \\ 
$\ 2:cab\ \ \ \ \ \mathbf{a}$ &  & $\ 2:cba\ \ \ \ \ $\textbf{not} $\mathbf{c%
}$ &  & $\ 2:bca\ \ \ \ \ $\textbf{not} $\mathbf{b}$ \\ 
$\ 3:bca$ &  & $\ 3:acb$ \ \ \ \ $not$ $a;$ $\ \therefore b$ &  & $\ 3:abc$
\\ 
&  &  &  &  \\ 
$\mathbf{3}$ &  & $\mathbf{7}$ &  & $\mathbf{11}$ \\ 
&  &  &  &  \\ 
$\ 1:acb$ &  & $\ 1:bca$ &  & $\ 1:cba$ \\ 
$\ 2:cba\ \ \ \ \ \mathbf{a}$ &  & $\ 2:cab$ \ \ \ \ $not$ $b$ &  & $\ 2:bac$
\\ 
$\ 3:bac$ &  & $\ 3:abc$ &  & $\ 3:acb$ \\ 
&  &  &  &  \\ 
$\mathbf{4}$ &  & $\mathbf{8}$ &  & $\mathbf{12}$ \\ 
&  &  &  &  \\ 
$\ 1:acb$ &  & $\ 1:bca$ &  & $\ 1:cba$ \\ 
$\ 2:bac\ \ \ \ \ \mathbf{a}$ &  & $\ 2:abc$ \ \ \ \ $not$ $b$ &  & $\ 2:acb$
\ \ \ \ $not$ $b$ \\ 
$\ 3:cba$ &  & $\ 3:cab$ &  & $\ 3:bac$%
\end{tabular}%
\]%
If $a$ were also chosen at profile 6, everyone would be decisive for $a$
against $b$ and $b$ would never be chosen on $NP(3,3)$. So $a$ is not chosen
there which means that $b$ is chosen there. But then $bD_{1}c$ and so $c$ is
not chosen at profiles 7 and 8. This means $a$ is chosen at each of those
profiles and that implies everyone is decisive for $a$ against $c$ and $c$
would never be chosen on $NP(3,3)$. Therefore our assumption that $a$ is
chosen at profile 5 fails and $b$ must be chosen there. Since $b$ is chosen
at profile 5, $bD_{1}c$ and thus $c$ is not chosen at profiles 7 and 8.
\pagebreak

\bigskip

\textbf{Profile 6}. 
\[
\begin{tabular}{lllll}
\underline{$\mathbf{1}$} &  & $\mathbf{5}$ &  & $\mathbf{9}$ \\ 
&  &  &  &  \\ 
$\ 1:abc$ &  & $\ 1:bac$ &  & $\ 1:cab$ \\ 
$\ 2:bca\ \ \ \ \ \mathbf{a}$ &  & $\ 2:acb\ \ \ \ \ \mathbf{b}$ &  & $\
2:abc$ \ \ \ \ \textbf{not} $\mathbf{b}$ \\ 
$\ 3:cab$ &  & $\ 3:cba$ &  & $\ 3:bca$ \\ 
&  &  &  &  \\ 
$\mathbf{2}$ &  & $\mathbf{6}$ &  & $\mathbf{10}$ \\ 
&  &  &  &  \\ 
$\ 1:abc$ &  & $\ 1:bac$ &  & $\ 1:cab$ \\ 
$\ 2:cab\ \ \ \ \ \mathbf{a}$ &  & $\ 2:cba\ \ \ \ \ $\textbf{not} $\mathbf{c%
};$ \ $a?$ &  & $\ 2:bca\ \ \ \ \ $\textbf{not} $\mathbf{b}$ \\ 
$\ 3:bca$ &  & $\ 3:acb$ &  & $\ 3:abc$ \\ 
&  &  &  &  \\ 
$\mathbf{3}$ &  & $\mathbf{7}$ &  & $\mathbf{11}$ \\ 
&  &  &  &  \\ 
$\ 1:acb$ &  & $\ 1:bca$ &  & $\ 1:cba$ \\ 
$\ 2:cba\ \ \ \ \ \mathbf{a}$ &  & $\ 2:cab$ \ \ \ \ \textbf{not} $\mathbf{c}%
;$ &  & $\ 2:bac$ \\ 
$\ 3:bac$ &  & $\ 3:abc$ \ \ \ \ $not$ $b;$ \ $\therefore a$ &  & $\ 3:acb$
\\ 
&  &  &  &  \\ 
$\mathbf{4}$ &  & $\mathbf{8}$ &  & $\mathbf{12}$ \\ 
&  &  &  &  \\ 
$\ 1:acb$ &  & $\ 1:bca$ &  & $\ 1:cba$ \\ 
$\ 2:bac\ \ \ \ \ \mathbf{a}$ &  & $\ 2:abc$ \ \ \ \ \textbf{not} $\mathbf{c}
$ &  & $\ 2:acb$ \\ 
$\ 3:cba$ &  & $\ 3:cab$ \ \ \ \ $not$ $b;$ \ $\therefore a$ &  & $\ 3:bac$%
\end{tabular}%
\]%
\medskip Suppose $a$ is chosen at profile 6, so that $aD_{3}b$. This implies 
$b$ is not chosen at profiles 7 and 8, so $a$ is chosen there. But then
everyone is decisive for $a$ against $c$ and $c$ would never be chosen on $%
NP(3,3)$. Accordingly, $a$ must not be chosen at profile 6, so $b$ must be.
Because $b$ is chosen at profiles 5 and 6, neither individual 2 nor 3 can be
decisive for $a$ against $b$, and thus $a$ cannot be selected at profiles 11
or 12. \pagebreak

\bigskip

\textbf{Profiles 7 and 8}.\medskip\ 
\[
\begin{tabular}{lllll}
$\mathbf{1}$ &  & $\mathbf{5}$ &  & $\mathbf{9}$ \\ 
&  &  &  &  \\ 
$\ 1:abc$ &  & $\ 1:bac$ &  & $\ 1:cab$ \\ 
$\ 2:bca\ \ \ \ \ \mathbf{a}$ &  & $\ 2:acb\ \ \ \ \ \mathbf{b}$ &  & $\
2:abc$ \ \ \ \ \textbf{not} $\mathbf{b}$ \\ 
$\ 3:cab$ &  & $\ 3:cba$ &  & $\ 3:bca$ \\ 
&  &  &  &  \\ 
$\mathbf{2}$ &  & $\mathbf{6}$ &  & $\mathbf{10}$ \\ 
&  &  &  &  \\ 
$\ 1:abc$ &  & $\ 1:bac$ &  & $\ 1:cab$ \\ 
$\ 2:cab\ \ \ \ \ \mathbf{a}$ &  & $\ 2:cba\ \ \ \ \ \mathbf{b}$ &  & $\
2:bca\ \ \ \ \ $\textbf{not} $\mathbf{b}$ \\ 
$\ 3:bca$ &  & $\ 3:acb$ &  & $\ 3:abc$ \ \ \ \ $not$ $c;$ \ $\therefore a$
\\ 
&  &  &  &  \\ 
$\mathbf{3}$ &  & $\mathbf{7}$ &  & $\mathbf{11}$ \\ 
&  &  &  &  \\ 
$\ 1:acb$ &  & $\ 1:bca$ &  & $\ 1:cba$ \\ 
$\ 2:cba\ \ \ \ \ \mathbf{a}$ &  & $\ 2:cab$ \ \ \ \ \textbf{not} $\mathbf{c}%
;$ \ $a?$ &  & $\ 2:bac\ \ \ \ \ $\textbf{not} $\mathbf{a}$ \\ 
$\ 3:bac$ &  & $\ 3:abc$ &  & $\ 3:acb$ \ \ \ \ $not$ $c;$ \ $\therefore b$
\\ 
&  &  &  &  \\ 
$\mathbf{4}$ &  & $\mathbf{8}$ &  & $\mathbf{12}$ \\ 
&  &  &  &  \\ 
$\ 1:acb$ &  & $\ 1:bca$ &  & $\ 1:cba$ \\ 
$\ 2:bac\ \ \ \ \ \mathbf{a}$ &  & $\ 2:abc$ \ \ \ \ \textbf{not} $\mathbf{c}
$ &  & $\ 2:acb\ \ \ \ \ $\textbf{not }$a$ \\ 
$\ 3:cba$ &  & $\ 3:cab$ &  & $\ 3:bac$ \ \ \ \ $not$ $c;$ \ $\therefore b$%
\end{tabular}%
\]%
Now suppose $a$ is chosen at profile 7 so that, by Lemma 4, $aD_{3}c$. Then $%
c$ is not chosen at profiles 10, 11, and 12. But $b$ being chosen at 11 and
12 implies $bD_{2}c$ and $bD_{3}c$, which in turn imply $c$ is never chosen
on $NP(3,3)$. So $a$ is not chosen at profile 7, so $b$ must be. Thus $%
bD_{1}a$ and $a$ is not chosen at profile 8, so $b$ is chosen
there.\pagebreak

\bigskip

\textbf{Profile 9}.\medskip

\qquad Suppose $a$ is chosen at profile 9. Then $aD_{2}c$ and so $c$ is also
not chosen at profiles 11 or 12 so $b$ is chosen there. But then $bD_{2}c$
and $bD_{3}c$. Thus $c$ is not in the range of $g$ restricted to $NP(3,3)$,
a contradiction. So we can't have $a$ chosen at profile 9 and thus $c$ is
chosen there.\medskip 
\[
\begin{tabular}{lllll}
$\mathbf{1}$ &  & $\mathbf{5}$ &  & $\mathbf{9}$ \\ 
&  &  &  &  \\ 
$\ 1:abc$ &  & $\ 1:bac$ &  & $\ 1:cab$ \\ 
$\ 2:bca\ \ \ \ \ \mathbf{a}$ &  & $\ 2:acb\ \ \ \ \ \mathbf{b}$ &  & $\
2:abc$ \ \ \ \ \textbf{not} $\mathbf{b};$ \ $a?$ \\ 
$\ 3:cab$ &  & $\ 3:cba$ &  & $\ 3:bca$ \\ 
&  &  &  &  \\ 
$\mathbf{2}$ &  & $\mathbf{6}$ &  & $\mathbf{10}$ \\ 
&  &  &  &  \\ 
$\ 1:abc$ &  & $\ 1:bac$ &  & $\ 1:cab$ \\ 
$\ 2:cab\ \ \ \ \ \mathbf{a}$ &  & $\ 2:cba\ \ \ \ \ \mathbf{b}$ &  & $\
2:bca\ \ \ \ \ $\textbf{not} $\mathbf{b}$ \\ 
$\ 3:bca$ &  & $\ 3:acb$ &  & $\ 3:abc$ \ \ \ \ $not$ $c;$ \ $\therefore a$
\\ 
&  &  &  &  \\ 
$\mathbf{3}$ &  & $\mathbf{7}$ &  & $\mathbf{11}$ \\ 
&  &  &  &  \\ 
$\ 1:acb$ &  & $\ 1:bca$ &  & $\ 1:cba$ \\ 
$\ 2:cba\ \ \ \ \ \mathbf{a}$ &  & $\ 2:cab$ \ \ \ \ $\mathbf{b}$ &  & $\
2:bac\ \ \ \ \ $\textbf{not} $\mathbf{a}$ \\ 
$\ 3:bac$ &  & $\ 3:abc$ &  & $\ 3:acb$ \ \ \ \ $not$ $c;$ \\ 
&  &  &  &  \\ 
$\mathbf{4}$ &  & $\mathbf{8}$ &  & $\mathbf{12}$ \\ 
&  &  &  &  \\ 
$\ 1:acb$ &  & $\ 1:bca$ &  & $\ 1:cba$ \\ 
$\ 2:bac\ \ \ \ \ \mathbf{a}$ &  & $\ 2:abc$ \ \ \ \ $\mathbf{b}$ &  & $\
2:acb\ \ \ \ \ $\textbf{not} $\mathbf{a}$ \\ 
$\ 3:cba$ &  & $\ 3:cab$ &  & $\ 3:bac$ \ \ \ \ $not$ $c;$%
\end{tabular}%
\]

\pagebreak

\bigskip

\textbf{Profiles 10, 11, and 12}.\medskip

\qquad Having $c$ chosen at profile 9 implies $cD_{1}b$, so $b$ is not
chosen at profiles 11 or 12, so $c$ is also chosen there. And if $a$ were
chosen at 10, we would have $aD_{3}c$, contrary to the choice of $c$ at
profile 11. So $c$ is also chosen at profile 11 and $\#1$ is a dictator on $%
VP$.\medskip

\qquad \textbf{Step 3}, the last stage in the proof of Theorem 3, is showing
that if strategy-proof $g$ on $NP(3,3)$ has range $X$, then so does the
restriction of $g$ to $VP$. Equivalently, we show that if the range of $g$
restricted to $VP$ has range of just two alternatives, say $\{a,b\}$, then $%
c $ is not in the range of $g$.\medskip

\qquad \textbf{Case A}. Suppose $a$ is chosen at some profile $u$ in $VP$
where only one person prefers $a$ to $b$ but $c$ is not chosen in $VP$.
Apart from permutations of individuals, $u$ must be\medskip 
\[
u:%
\begin{tabular}{|c|c|c|}
\hline
1 & 2 & 3 \\ \hline
$a$ & $c$ & $b$ \\ 
$c$ & $b$ & $a$ \\ 
$b$ & $a$ & $c$ \\ \hline
\end{tabular}%
\]

Thus, $aD_{1}b$ and so $a$ is chosen at the 6 out of 12 profiles in $VP$
where $\#1$ prefers $a$ to $b$.\medskip 
\[
\begin{tabular}{lllll}
$\mathbf{1}$ &  & $\mathbf{5}$ &  & $\mathbf{9}$ \\ 
$\ 1:abc$ &  & $\ 1:bac$ &  & $\ 1:cab$ \\ 
$\ 2:bca\ \ \ \ \ a$ &  & $\ 2:acb$ \ \ \ \  &  & $\ 2:abc$ \ \ \ \ $a$ \\ 
$\ 3:cab$ &  & $\ 3:cba$ &  & $\ 3:bca$ \\ 
&  &  &  &  \\ 
$\mathbf{2}$ &  & $\mathbf{6}$ &  & $\mathbf{10}$ \\ 
$\ 1:abc$ &  & $\ 1:bac$ &  & $\ 1:cab$ \\ 
$\ 2:cab\ \ \ \ \ a$ &  & $\ 2:cba$ \ \ \ \  &  & $\ 2:bca\ \ \ \ \ a$ \\ 
$\ 3:bca$ &  & $\ 3:acb$ &  & $\ 3:abc$ \\ 
&  &  &  &  \\ 
$\mathbf{3}$ &  & $\mathbf{7}$ &  & $\mathbf{11}$ \\ 
$\ 1:acb$ &  & $\ 1:bca$ &  & $\ 1:cba$ \\ 
$\ 2:cba\ \ \ \ \ a$ &  & $\ 2:cab$ \ \ \ \  &  & $\ 2:bac$ \\ 
$\ 3:bac$ &  & $\ 3:abc$ &  & $\ 3:acb$ \\ 
&  &  &  &  \\ 
$\mathbf{4}$ &  & $\mathbf{8}$ &  & $\mathbf{12}$ \\ 
$\ 1:acb$ &  & $\ 1:bca$ &  & $\ 1:cba$ \\ 
$\ 2:bac\ \ \ \ \ a$ &  & $\ 2:abc$ \ \ \ \  &  & $\ 2:acb$ \\ 
$\ 3:cba$ &  & $\ 3:cab$ &  & $\ 3:bac$%
\end{tabular}%
\]%
\medskip By profile 2, $aD_{1}c$, by profile 9, $aD_{2}c$ and by profile 10, 
$aD_{3}c$. Lemma 4-4 implies we must have $c$ is not in the range of $g$%
.\pagebreak

\bigskip

\qquad \textbf{Case B}. Alternative $a$ is never chosen when only one person
prefers $a$ to $b$ under our supposition that $c$ does not belong to the
range of $g$ restricted to $VP$.\medskip 
\[
\begin{tabular}{lllll}
$\mathbf{1}$ &  & $\mathbf{5}$ &  & $\mathbf{9}$ \\ 
&  &  &  &  \\ 
$\ 1:abc$ &  & $\ 1:bac$ &  & $\ 1:cab$ \\ 
$\ 2:bca$ &  & $\ 2:acb$ \ \ \ \ $b$ &  & $\ 2:abc$ \\ 
$\ 3:cab$ &  & $\ 3:cba$ &  & $\ 3:bca$ \\ 
&  &  &  &  \\ 
$\mathbf{2}$ &  & $\mathbf{6}$ &  & $\mathbf{10}$ \\ 
&  &  &  &  \\ 
$\ 1:abc$ &  & $\ 1:bac$ &  & $\ 1:cab$ \\ 
$\ 2:cab$ &  & $\ 2:cba$ \ \ \ \ $b$ &  & $\ 2:bca$ \\ 
$\ 3:bca$ &  & $\ 3:acb$ &  & $\ 3:abc$ \\ 
&  &  &  &  \\ 
$\mathbf{3}$ &  & $\mathbf{7}$ &  & $\mathbf{11}$ \\ 
&  &  &  &  \\ 
$\ 1:acb$ &  & $\ 1:bca$ &  & $\ 1:cba$ \\ 
$\ 2:cba\ \ \ \ \ b$ &  & $\ 2:cab$ &  & $\ 2:bac$ \ \ \ \ $b$ \\ 
$\ 3:bac$ &  & $\ 3:abc$ &  & $\ 3:acb$ \\ 
&  &  &  &  \\ 
$\mathbf{4}$ &  & $\mathbf{8}$ &  & $\mathbf{12}$ \\ 
&  &  &  &  \\ 
$\ 1:acb$ &  & $\ 1:bca$ &  & $\ 1:cba$ \\ 
$\ 2:bac\ \ \ \ \ b$ &  & $\ 2:abc$ &  & $\ 2:acb$ \ \ \ \ $b$ \\ 
$\ 3:cba$ &  & $\ 3:cab$ &  & $\ 3:bac$%
\end{tabular}%
\]%
By profile 5, $bD_{1}c$; by profile 4, $bD_{2}c$; and by profile 3, $bD_{3}c$%
. Then Lemma 4 implies $c$ is not in the range of $g$. \ \ \ \ $\square $%
\bigskip

\section{Induction on n.\protect\medskip}

\qquad We have established $SP(n,m)$ for $n=3$ and $m=3$. We want to extend
this to more individuals (but still for $m=3$) by induction on $n$. The
analysis of this section calls for defining rules on $NP(n,3)$ as
retrictions of rules on $NP(n+1,3)$ to subdomains. In particular, we want to
define two rules on $NP(n,3)$ by considering their images at a profile $%
p=(p(1),...,p(n))$:\medskip

\qquad \qquad \qquad 
\begin{tabular}{l}
$g^{\ast }(p)=g(p(1),p(2),...,p(n-1),p(n),p(n))$ and\medskip \\ 
$g^{\ast \ast }(p)=g(p(1),p(1),p(2),...,p(n-1),p(n))$\medskip%
\end{tabular}

We will need these rules to be of full range when $g$ is of full range.
Accordingly, we define two subdomains of $NP(n+1,3)$:\medskip

\qquad $NP^{\ast }(n+1,3)$ is the subdomain of $NP(n+1,3)$ consisting of all 
$p\in NP(n+1,3)$ such that $p(n)=p(n+1)$;\medskip

\qquad $NP^{\ast \ast }(n+1,3)$ is the subdomain of $NP(n+1,3)$ consisting
of all $p\in NP(n+1,3)$ such that $p(1)=p(2)$.\medskip

Recall that $SP(n,3)$ is the statement: \textquotedblleft \textbf{every
strategy-proof rule on }$NP(n,3)$\textbf{\ is dictatorial if it has full
range}.\textquotedblright \medskip

\qquad \textbf{Theorem 6-1.} $SP(n,3)$ implies $SP(n+1,3)$.\medskip

\qquad \textbf{Proof:} Assume that $SP(n,3)$ is true. Let g be a given
strategy-proof social choice function on $NP(n+1,3)$ that has full
range.\medskip

\qquad We need to carry out several steps:\medskip

\qquad \textbf{1.} Given the set $N=\{1,2,...,n,n+1\}$ of individuals, we let

$\qquad \qquad N^{\ast }=\{1,2,...,n\}$.\medskip

\qquad \textbf{2.} Associate with each profile $p$ on $X$ in $NP(n,3)$ two
profiles,

$\qquad \qquad p^{\ast }$ and $p^{\ast \ast }$ on $X$ in $NP(n+1,3)$.\medskip

\qquad \textbf{3.} Define $g^{\ast }$ on $NP(n,3)$ by relating $g^{\ast }(p)$
to $g(p^{\ast })$ and define

$\qquad \qquad g^{\ast \ast }$ on $NP(n,3)$ by relating $g^{\ast \ast }(p)$
to $g(p^{\ast \ast })$.\medskip

\qquad \textbf{4.} Show that strategy-proofness of $g$ on $NP(n+1,3)$ implies

\qquad \qquad strategy-proofness of $g^{\ast }$ and $g^{\ast \ast }$ on $%
NP(n,3)$.\medskip

\qquad \textbf{5.} Show that Range($g$) $=X$ implies Range($g^{\ast }$) $=X$
and

\qquad \qquad Range($g^{\ast \ast }$) $=X$. This step is carried out in
companion paper

\qquad \qquad Campbell and Kelly, 2014b).\medskip

\qquad \textbf{6.} Conclude by the induction hypothesis that $g^{\ast }$ and 
$g^{\ast \ast }$ are dictatorial.\medskip

\qquad \textbf{7.} Show that the dictatorship of $g^{\ast }$ and $g^{\ast
\ast }$ implies dictatorship for $g$.\medskip

For the second step, given a profile $p=(p(1),p(2),...,p(n))$ in $NP(n,3)$,
let\medskip

\[
p^{\ast }=(p(1),p(2),...,p(n),p(n)) 
\]

where the last two orderings are the same, and

\[
p^{\ast \ast }=(p(1),p(1),p(2),...,p(n)) 
\]

where the first two orderings are the same.\medskip

\qquad \qquad \textbf{Remark}: $p$ is in NP(n, 3) if and only if $p^{\ast }$
is in $NP(n+1,3)$; and $p$ is in $NP(n,3)$ if and only if $p^{\ast \ast }$
is in $NP(n+1,3)$.\medskip

\qquad For the third step, we define $g^{\ast }$ and $g^{\ast \ast }$ as at
the beginning of this section.\medskip

\qquad For the fourth step, we need to show $g^{\ast }$ and $g^{\ast \ast }$
are strategy-proof. Consider rule $g\ast $: If individual $i<n$ can
manipulate $g\ast $ at some profile $(p(1),p(2),...,p(n))$ then that
individual can obviously manipulate g at profile $(p(1),p(2),...,p(n),p(n))$
because \medskip

\qquad \qquad $g^{\ast }(p(1),p(2),...,p(n))=g(p(1),p(2),...,p(n),p(n))$%
.\medskip

Next, we show that the remaining individual, $n$, cannot manipulate $g^{\ast
}$. Suppose to the contrary that $n$ manipulates $g^{\ast }$ at profile $%
p=(p(1),p(2),...,p(n-1),p(n))$ in $NP(n,3)$ with $g^{\ast }(p)=x$, by
submitting ordering $v$ and getting profile $q=(p(1),p(2),...,p(n-1),v)$ in $%
NP(n,3)$ with $g^{\ast }(q)=y\neq x$ and $y\succ _{p(n)}x$.\medskip

\qquad \qquad \qquad $g(p(1),p(2),...,p(n-1),p(n),p(n))=x$; and

\qquad \qquad \qquad $g(p(1),p(2),...,p(n-1),v,v)=y$. \medskip

Note by the remark above, that these profiles are in $NP(n+1,3)$.\medskip

Construct profile $w=(p(1),p(2),...,p(n-1),p(n),v)$. Note that $w$ is in $%
NP(n+1,3)$. Now $x\succ _{p(n)}g(w)$ or $g$ is manipulable by $n+1$ at $%
(p(1),p(2),...,p(n-1),p(n),p(n))$, and $g(w)\succ _{p(n)}y$ or g is
manipulable by $n$ at w. Transitivity then implies $x\succ _{p(n)}y$,
contrary to our earlier assumption of manipulability of $g^{\ast }$.
Therefore, $g^{\ast }$ is strategy-proof.\medskip

\qquad An analogous argument establishes the strategy-proofness of $g^{\ast
\ast }$.\medskip

\qquad For the fifth step, as noted above, Range($g^{\ast }$) = Range($%
g^{\ast \ast }$) $=X$ is proven as the \textbf{n-Range Lemma} in (Campbell
and Kelly, 2014b).\medskip

\qquad It follows from Steps 4 and 5, and our induction hypothesis, that $%
g^{\ast }$ and $g^{\ast \ast }$ are dictatorial. We now want to prove that
existence of a dictator for each of $g^{\ast }$ and for $g^{\ast \ast }$
enables us to determine a dictator for $g$.\medskip

\qquad \textbf{Case 1}. The dictator for $g^{\ast }$ is an individual $j<n$.
Without loss of generality, $j=1$. Let $x$ be an arbitrary member of $X$,
and let

\[
q=(q(1),q(2),p(3),...,q(n),q(n+1)) 
\]%
be any profile in $NP(n+1,3)$ with the topmost element of q(1) being $x$.
Let $p=(q(1),p(2),...,p(n-1),p(n),p(n+1))$ be a profile with $p(1)=q(1)$ and
for all $i>1$, $p(i)$ is the inverse of $p(1$). Note that $p\in NP(n,3)$.
Now $g(p)=x$ since $p(n)=p(n+1)$ and $\#1$ is a dictator for $g^{\ast }$%
.\medskip

Next we consider a standard sequence from $p$ to $q$:%
\begin{eqnarray*}
p^{1} &=&p=(p(1),p(2),p(3),...,p(n),p(n))=(q(1),p(2),...,p(n),p(n))\medskip
\\
p^{2} &=&(q(1),q(2),p(3),...,p(n),p(n))\medskip \\
&&\vdots \\
p^{n} &=&(q(1),q(2),q(3),...,q(n),p(n))\medskip \\
p^{n+1} &=&q=(q(1),q(2),q(3),...,q(n),q(n+1)).\medskip
\end{eqnarray*}

Note that every profile in this sequence is an element of $NP(n+1,3)$, and
that $g(p)=x$. Then also $g(p^{2})=x$ since otherwise, with $x$ at the
bottom of $p(2)$, $g$ would be manipulable by $\#2$ at $p^{1}$. Similarly, $%
g(p^{3})=x$ or $g$ would be manipulable by $\#3$ at $p^{2}$. Continuing in
this fashion, $g(q)=x$ and so $\#1$ is a dictator for $g$.\medskip

\qquad \textbf{Case 2}. The dictator for $g^{\ast \ast }$ is an individual $%
j>1$. Using the argument of Case 1 as a template, we can show that $j+1$ is
a dictator for $g$.\medskip

\qquad \textbf{Case 3}. The dictator for $g^{\ast }$ is individual $n$ and
the dictator for $g^{\ast \ast }$ is $\#1$. Consider a profile $u$ such that 
$u(1)=u(2)=xy...z$ and $u(n)=u(n+1)=z...yx$. This profile is in $NP(n+1,3)$.
Since $\#1$ is a dictator for $g^{\ast \ast }$, $g(u)=x$; since $n$ is a
dictator for $g^{\ast }$, $g(u)=z$. This contradiction shows Case 3 isn't
possible. \ \ \ \ $\square $\medskip

\section{Induction on m.\protect\medskip}

We have established that $SP(n,3)$ holds for all $n\geq 3$. We now show%
\textbf{\medskip }

\textbf{Theorem 6-1}. \ For arbitrary $n\geq 3$ and arbitrary $m\geq 3$ the
statement $SP(n,m)$ implies $SP(n,m+1)$. \textbf{\medskip }

\textbf{Proof}: \ Assume that $SP(n,m)$ is true.\textbf{\medskip }

\qquad Let $g$ be a given strategy-proof social choice function on $%
NP(n,m+1) $ that has full range. Now we define a rule $g^{\ast }$ based on $%
g $. Select arbitrary, but distinct, $w$ and $z$ in $X$. Let $NP^{wz}(n,m+1)$
be the set of profiles in $NP(n,m+1)$ such that alternatives $w$ and $z$ are
contiguous in each individual ordering. Choose some alternative $x^{\ast }$
that does not belong to $X$ and set $X^{\ast }=\{x^{\ast }\}\cup X\backslash
\{w,z\}$. Then $g^{\ast }$ will have domain $D^{\ast }$ by which we mean the
domain $NP(n,m)$ when the feasible set is $X^{\ast }$. To define $g^{\ast }$
we begin by selecting arbitrary profile $p\in D^{\ast }$, and then we choose
some profile $r\in NP^{wz}(n,m+1)$ such that\textbf{\medskip }

\qquad 1.\qquad $r|X\backslash \{w,z\}=p|X\backslash \{w,z\}$, and

\qquad 2.\qquad for any $i\in \{1,2,...,n\}$, we have 
\[
\{x\in X\backslash \{w,z\}:x\succ _{r(i)}w\}=\{x\in X\backslash
\{w,z\}:x\succ _{p(i)}x^{\ast }\}. 
\]

In words, we create $r$ from $p$ by replacing $x^{\ast }$ with $w$ and $z$
so that $w$ and $z$ are contiguous in each $r(i)$, and $r$ does not exhibit
any Pareto domination, and in each $r(i)$ either $w$ or $z$ occupies the
same rank as $x^{\ast }$ in $p(i)$. We next show that the selected
alternative, which we can denote $f(p)$, is independent of the choice of
profile r.

Suppose that $s$ and $t$ both belong to $NP^{wz}(n,m+1)$ and both 1 and 2
hold for $r=s$ and $r=t$. We will show that if $t\neq s$ there is a profile $%
r^{\prime }$ in $NP^{wz}(n,m+1)$ and an individual $h$ such that

\qquad \qquad $r^{\prime }(i)=s(i)$ for all $i\neq h$ and $r^{\prime
}(h)=t(h)\neq s(h)$ and

\qquad \qquad $g(r^{\prime })=g(s)$ if $g(s)\in X\backslash \{w,z\}$, and $%
g(r^{\prime })\in \{w,z\}$ if $g(s)\in \{w,z\}$.

That will establish that we can create $t$ from $s$ in a serious of stages
without changing the selected alternative unless it changes from $w$ to $z$
or from $z$ to $w$. \ We have shown that $g^{\ast }(p)$ is independent of
the choice of $r$ in $NP^{wz}(n,m+1)$, provided that 1 and 2 are both
satisfied. Thus $g^{\ast }$ is well defined.

\qquad Choose some $j$ in $\{1,2,...,n\}$ such that $s(j)\neq t(j)$. Without
loss of generality, assume that $w\succ _{t(j)}z$. Create $u$ from $s$ by
switching $w$ and $z$ in $s(j)$. If $u\in NP^{wz}(n,m+1)$ then the fact that 
$w$ and $z$ are contiguous in $s$ and strategy-proofness of $g$ imply that $%
g(u)=g(s)$ if $g(s)\in X\backslash \{w,z\}$, and $g(u)\in \{w,z\}$ if $%
g(s)\in \{w,z\}$. Set $r^{\prime }=u$.

\qquad Suppose that $u\notin NP(n,m+1)$. Then we must have $z\succ _{t(i)}w$
for all $i\neq j$. Because $z$ does not Pareto dominate $w$ at $t$ and we
have $t(j)\neq s(j)$, and $u(i)=s(i)$ for all $i\neq j$ there exists an
individual $h$ such that $z\succ _{s(h)}w$ and $t(h)\neq s(h)$. Then switch $%
w$ and $z$ in s(h) to form $v\in NP^{wz}(n,m+1)$. Strategy-proofness of $g$
implies that $g(v)=g(s)$ if $g(s)\in X\backslash \{w,z\}$, and $g(v)\in
\{w,z\}$ if $g(s)\in \{w,z\}$. Set $r^{\prime }=v$.\textbf{\medskip }

Therefore, we can use arbitrary $f$ satisfying 1 and 2 above to define $%
g^{\ast }$ as follows.\textbf{\medskip }

\textbf{Definition}: $g^{\ast }(p)=g(f(p))$ if $g(f(p))\in X\backslash
\{w,z\}$ and $g^{\ast }(p)=x^{\ast }$ if $g(f(p))\in \{w,z\}$.\textbf{%
\medskip }

The rule $g^{\ast }$ is obviously strategy-proof if $g$ is. The \textbf{%
m-Range Lemma} in the companion paper (Campbell and Kelly, 2014b)
establishes that the range of $g^{\ast }$ is $X^{\ast }$. We can now claim
that $g^{\ast }$ is dictatorial using the induction hypothesis.\textbf{%
\medskip }

\qquad Without loss of generality, assume that person $1$ is the dictator
for $g^{\ast }$. We will show that person 1 is a dictator for $g$ and that
will establish $SP(n,m+1)$. We prove dictatorship of $g$ by showing that,
for arbitrary $x\in X$, $g(p)=x$ for any profile in $NP(n,m+1)$ that has $x$
at the top of person 1's preference ordering. The first step is to show that
if $x$ is selected at some profile that has $x$ at the top of person 1's
ordering and at the bottom of everyone else's, then $g$ will select $x$ at
any profile for which it is the top ranked alternative of person 1. Again,
the analysis is complicated by the need to stay inside $NP(n,m+1)$.\textbf{%
\medskip }

\textbf{Lemma 6-2}: For any $x\in X$, if there is a profile $r\in NP(n,m)$
such that $g(r)=x$, and $x$ is at the top of $r(1)$ but at the bottom of $%
r(i)$ for all $i>1$ then $g(p)=x$ for any profile $p$ in $NP(n,m+1)$ such
that $x$ is at the top of $p(1)$.\textbf{\medskip }

\textbf{Proof}: Let $p$ be an arbitrary member of $NP(n,m+1)$ such that $x$
is at the top of $p(1)$. By hypothesis, there exists a profile $r\in
NP(n,m+1)$ such that $g(r)=x$, and $x$ is at the top of $r(1)$ and the
bottom of $r(i)$ for all $i>1$. This justifies the first in the following
sequence of statements, and each of the others follows from its predecessor
by strategy-proofness of $g$, as we explain at the end of the list.\textbf{%
\medskip }

\qquad 1.$\qquad g(r(1),r(2),r(3),...,r(n-1),r(n))=x\mathbf{\medskip }$

\qquad 2.$\qquad g(r(1),r(1)^{-1},r(1)^{-1},...,r(1)^{-1},r(1)^{-1})=x%
\mathbf{\medskip }$

\qquad 3.\qquad $g(r(1),p(1)^{-1},p(1)^{-1},...,p(1)^{-1},r(1)^{-1})=x$%
\textbf{\medskip }

\qquad 4.\qquad $g(p(1),p(1)^{-1},p(1)^{-1},...,p(1)^{-1},r(1)^{-1})=x$%
\textbf{\medskip }

\qquad 5.\qquad $g(p(1),p(1)^{-1},p(3),...,p(n-1),p(n))=x$\textbf{\medskip }

\qquad 6.\qquad $g(p(1),p(2),p(3),...,p(n-1),p(n))=x$.\textbf{\medskip }

Note that each of the six profiles belongs to $NP(n,m+1)$. Statement 1
implies 2 because we can replace $r(2)$ with $r(1)^{-1}$. Alternative $x$
will still be selected because $x$ is selected at the previous stage and $x$
is bottom ranked by $r(2)$. Then replace $r(3)$ with $r(1)^{-1}$.
Alternative $x$ will still be selected because $x$ is selected at the
previous stage and $x$ is bottom ranked by $r(3)$. And so on. Statement 2
implies 3 because we can replace person $i$'s ordering $r(i)^{-1}$ with $%
p(1)^{-1}$ one individual at a time, for $2\leq i\leq n-1$. Because $x$ is
bottom ranked by $r(i)^{-1}$, if $x$ is not selected after replacing $%
r(i)^{-1}$ with $p(1)^{-1}$ for individual $i$ that person could manipulate
at the profile for which $r(1)^{-1}$ is his true preference. Statement 3
implies 4 because the fourth profile is obtained from the third by replacing 
$r(1)$ with $p(1)$, and $x$ is top ranked by $p(1)$ so person 1 could
manipulate at the fourth profile if $x$ were not selected there. Statement 4
implies 5 because we can replace person $i$'s ordering with $p(i)$ one
individual at a time, for $i\geq 3$. Because $x$ is bottom ranked by $%
p(1)^{-1}$ and $r(1)^{-1}$, if $x$ is not selected after replacing $i$'s
ordering at profile 4 that person could manipulate at the pre-replacement
profile. Finally, statement 5 implies 6 because the last profile is obtained
from its predecessor by changing person 2's ordering, replacing $p(1)^{-1}$
with $p(2)$ and $x$ is at the bottom of $p(1)^{-1}$. \ \ \ \ $\square $%
\textbf{\medskip }

\qquad Now we establish that there is a profile $r$ in $NP(n,m+1)$ such that 
$g(r)=x$ and $x$ is at the top of $r(1)$ and at the bottom of $r(i)$ for all 
$i>1$.\textbf{\medskip }

\textbf{Lemma 6-3}: For each $x$ $\in $ $X$ there is a profile $p\in
NP(n,m+1)$ such that $g(p)=x$, and $x$ is at the top of $p(1)$ but at the
bottom of $p(i)$ for all $i>1$.\textbf{\medskip }

\textbf{Proof}: \textbf{Step 1}: We show that, for arbitrary $x\in
X\backslash \{w,z\}$, we have $g(p)=x$ for \textit{any} $p\in NP(n,m+1)$
such that $x$ is at the top of $p(1)$.\textbf{\medskip }

\qquad Let $\alpha $ be any member of $L(X)$ such that $x$ is at the top of $%
\alpha $ and $w$ and $z$ are contiguous in $\alpha $. Set $r=(\alpha ,\alpha
^{-1},\alpha ^{-1},...,\alpha ^{-1})$. We have $g(r)=x$ because person 1
dictates $g^{\ast }$. By Lemma 6-2, $g(p)=x$ for any $p\in NP(n,m+1)$ such
that $x$ is at the top of $p(1)$.\textbf{\medskip }

\qquad The previous paragraph will not suffice when $x\in \{w,z\}$ because
the dictatorship of $g^{\ast }$ only establishes that $g$ selects either $w$
or $z$ when those two alternatives rank first and second (in either order)
in person 1's preference scheme. We need to show that $g$ selects $w$
(resp., z) when $w$ (resp., $z$) ranks at the top for person 1.\textbf{%
\medskip }

\textbf{Step 2}: We show that if person 1's ordering has a member of $%
\{w,z\} $ at the top then $g$ will select 1's first or second ranked
alternative.\textbf{\medskip }

\qquad Let $q$ be any profile in $NP(n,m+1)$ such that $z$ ranks first in $%
q(1)$ and $w$ ranks second, or $w$ ranks first in $q(1)$ and $z$ ranks
second. (Note that $q$ does not have to belong to $NP^{wz}(n,m+1)$.) \ We
show that $g(q)\in \{w,z\}$. Let $r$ be the profile for which $r(1)=q(1)$,
and $r(i)=q(1)^{-1}$ for all $i\geq 2$. We have $g(r)\in \{w,z\}$ because
person 1 dictates $g^{\ast }$. Then $g(r(1),q(2),q(3),...,q(n-1),r(n))\in
\{w,z\}$ by strategy-proofness of $g$. (For $2\leq i\leq n-1$ replace $r(i)$
with $q(i)$ one individual at a time. The fact that $r(n)=r(1)^{-1}$
guarantees that we remain in $NP(n,m+1)$.) We have $%
g(q)=g(r(1),q(2),q(3),...,q(n-1),q(n))\in \{w,z\}$ by strategy-proofness of $%
g$.\textbf{\medskip }

\qquad Now, let $p$ be any profile in $N(n,m+1)$ such that $p(1)=zx...$ for
some $x\in X\backslash \{w,z\}$. We prove $g(p)\in \{z,x\}$. Let $\beta $ be
any member of L(X) such that $x$ ranks first and z ranks second. Let r
denote the profile $(\beta ,\beta ^{-1},\beta ^{-1},...,\beta ^{-1})$. We
have $g(r)=x$ by Step 1. This is the first in the following sequence of
statements. Each of statements 2 through 5 follows from its predecessor by
strategy-proofness of $g$, as we explain in the paragraph following the list.%
\textbf{\medskip }

\qquad 1.\qquad $g(r(1),r(1)^{-1},r(1)^{-1},...,r(1)^{-1},r(1)^{-1})=x$%
\textbf{\medskip }

\qquad 2.\qquad $g(r(1),r(1)^{-1},p(1)^{-1},...,p(1)^{-1},p(1)^{-1})=x$%
\textbf{\medskip }

\qquad 3.\qquad $g(p(1),r(1)^{-1},p(1)^{-1},...,p(1)^{-1},p(1)^{-1})\in
\{x,z\}$\textbf{\medskip }

\qquad 4.\qquad $g(p(1),p(2),p(3),...,p(n-1),p(1)^{-1})\in \{x,z\}$\textbf{%
\medskip }

\qquad 5.\qquad $g(p(1),p(2),p(3),...,p(n-1),p(n))\in \{x,z\}$\textbf{%
\medskip }

Statement 1 implies 2 because, for $i\geq 3$, we can replace person $i$'s
ordering $r(i)=r(1)^{-1}$ with $p(1)^{-1}$ one individual at a time. Then $x$
will still be selected at each stage because it is at the bottom of $r(i)$,
and if it were not selected at the new profile then person $i$ could
manipulate at the previous profile. Statement 2 implies 3 by
strategy-proofness because we replace $r(1)$ with $p(1)$ which has $z$
ranked first and $x$ second. Statement 3 implies 4 because, for $2\leq i\leq
n-1$, we replace $r(i)=r(1)^{-1}$ or $p(1)^{-1}$ with $p(i)$, one individual
at a time. Then strategy-proofness of $g$ implies that $x$ or $z$ will still
be selected because $x$ ranks at the bottom of $r(i)=r(1)^{-1}$ with $z$
second last and $p(1)^{-1}$ has $z$ ranked last and $x$ second last.
Statement 4 implies 5 because we replaced person $n$'s ordering, $p(1)^{-1}$%
, which has $z$ ranked last and $x$ second last.\textbf{\medskip }

\qquad Similarly, if $p(1)=(w,x,...)$ then $g(p)\in \{w,x\}$.\textbf{%
\medskip }

\textbf{Step 3}: Suppose that there exists $t\in NP(n,m+1)$ such that $%
g(t)=z $ and $z$ is at the top of $t(1)$. We show that there exists a
profile $u\in NP(n,m+1)$ such that $g(u)=z$ and $z$ is at the top of $u(1)$
and at the bottom of $u(i)$ for all $i>1$. And if there exists $t^{\prime
}\in NP(n,m+1) $ such that $g(t^{\prime })=w$ and $w$ is at the top of $%
t^{\prime }(1)$ then there is a profile $u^{\prime }\in NP(n,m+1)$ such that 
$g(u^{\prime })=w$ and $w$ is at the top of $u^{\prime }(1)$ and at the
bottom of $u^{\prime }(i)$ for all $i>1$.\textbf{\medskip }

\qquad By hypothesis, there exists a profile $t\in NP(n,m+1)$ such that $%
g(t)=z$ and $t(1)=zxy...$. That is, $z$ ranks first in $t(1)$, and the
second ranking alternative is denoted by $x$, with $y$ third. Without loss
of generality we assume that $z\succ _{t(i)}x$ for all $i\leq k$ (possibly $%
k=1$) and $x\succ _{t(i)}z$ for all $i>k$. We have $k<n$ because $z$ does
not Pareto dominate $x$. Define $p\in NP(n,m+1)$ by setting $p(i)=t(i)$ for
all $i\leq k$ and $p(i)=t(1)^{-1}$ for all $i>k$. If we change $t(i)$ for $%
i>k$ one individual at a time, alternative $z$ will be selected at each
stage. (Step 2 implies that either $z$ or $x$ will be selected, and if $x$
is selected then the individual whose ordering has changed could manipulate
at the previous profile.) Therefore, $g(p)=z$. If $k=1$ set $u=p$.\textbf{%
\medskip }

\qquad If $k>1$ create $q$ from $p$ by setting $q(2)=y...zx$, with $%
q(i)=p(i) $ for all $i\neq 2$. (That is, $q(2)$ is any ordering with $y$ on
top, $x$ last, and $z$ second last.) We have $q\in NP(n,m+1)$ because $%
p(n)=t(1)^{-1}$. And $g(q)\in \{z,x\}$ by Step 2. If $g(q)=x$ then person 2
can manipulate at $q$ via $p(2)$. Thus, $g(q)=z$.\textbf{\medskip }

\qquad Recall that $q(n)=p(n)=t(1)^{-1}$. Create $r$ from $q$ by switching
the order of $y$ and $x$ in $q(n)$, with $r(i)=q(i)$ for all $i\neq n$. We
have $r(n)=$ $...xyz$ and $y\succ _{q(2)}x$. Therefore, $r\in NP(n,m+1)$.
Step 2 implies that $g(r)\in \{z,x\}$. If $g(r)=x$ then person $n$ can
manipulate at $q$ via $r(n)$. Therefore, $g(r)=z$. \textbf{\medskip }

\qquad Create $s$ from $r$ by switching the order of $y$ and $x$ in $%
r(1)=t(1)$, with $s(i)=r(i)$ for all $i>1$. We have $s\in NP(n,m+1)$ because 
$x\succ _{r(n)}y$. Strategy-proofness implies $g(s)=z$.\textbf{\medskip }

\qquad Now create $u^{2}$ from $s$ by switching the order of $z$ and $x$ in $%
s(2)$. We have $u^{2}(2)=y...xz$, and $u^{2}\in NP(n,m+1)$ because $z\succ
_{s(1)}x$. Step 2 implies $g(u^{2})\in \{z,y\}$, but if $g(u^{2})=y$ then
person 2 can manipulate at $s$ via $u^{2}(2)$. Therefore, $g(u^{2})=z$. If $%
k=2$ we have $z$ at the top of $u^{2}(1)$ and at the bottom of $u^{2}(i)$
for all $i>1$, in which case we can set $u=u^{2}$. The following table
specifies $t(1)$ for the given profile $t$ and also summarizes the changes
that have been made. It only identifies the ordering for the individuals
whose preferences have altered. We have established that alternative $z$ is
selected at each stage. \textbf{\medskip }

Sketch of proof that u exists when $k=2$%
\begin{eqnarray*}
t(1) &=&zxy...\mathbf{\medskip } \\
p(i) &=&t(1)^{-1}\text{ }for\text{ }all\text{ }i>k\mathbf{\medskip } \\
q(2) &=&y...zx\mathbf{\medskip } \\
r(n) &=&...xyz\mathbf{\medskip } \\
s(1) &=&zyx,...\mathbf{\medskip } \\
u^{2}(2) &=&y...xz\mathbf{\medskip }
\end{eqnarray*}

\qquad If $k>2$ create $u^{3}$ from $u^{2}$ by replacing $u^{2}(3)$ with $%
u^{2}(2)$. Strategy proofness implies that $g(u^{3})=z$. Continue in this
fashion until we have $g(u^{k})=z$, with $z$ on top of $u^{k}(1)$ and at the
bottom of $u^{k}(i)$ for all $i>1$. Then set $u=u^{k}$.\textbf{\medskip }

\qquad Similarly, we can prove the existence of a profile $u^{\prime }$ in $%
NP(n,m+1)$ such that $g(u^{\prime })=w$, with $w$ at the top of $u^{\prime
}(1)$ and at the bottom of $u^{\prime }(i)$ for all $i>1$.\textbf{\medskip }

\textbf{Step 4}: We show that, for arbitrary $a\in \{w,z\}$, there does in
fact exist a profile $t\in NP(n,m+1)$ such that $g(t)=a$ and $a$ is at the
top of $t(1)$. Without loss of generality, $a=z$.\textbf{\medskip }

\qquad Alternative $z$ is in the range of $g$ so there is a profile $q\in
NP(n,m+1)$ such that $g(q)=z$. Suppose that $z$ is not the top alternative
of $q(1)$. Step 1 implies that the top alternative of $q(1)$ does not belong
to $X\backslash \{w,z\}$. Therefore, $q(1)=(w,...)$. Step 2 implies that the
alternative ranked second by $q(1)$ must be $z$. Therefore, $q(1)=wz...$ and 
$g(q)=z$. If $w\succ _{q(j)}z$ for some $j>1$ then let $r\in NP(n,m+1)$ be
the profile for which $r(i)=q(i)$ for all $i>1$ and $r(1)$ is obtained from $%
q(1)$ by switching $w$ and $z$. Then $g(r)=z$ by strategy-proofness. Set $%
t=r $.\textbf{\medskip }

\qquad Suppose that $z\succ _{q(i)}w$ for all $i>1$. Create $p$ from $q$ by
moving $z$ down in $q(2)$ until it's just above $w$ in person 2's ordering,
preserving%
\[
p(2)|X\backslash \{z\}=q(2)|X\backslash \{z\}, 
\]%
with $p(i)=q(i)$ for all $i\neq 2$. (Then $p=q$ if $w$ and $z$ are
contiguous in $q(2)$.) We have $p\in NP(n,m+1)$ because $z\succ _{q(1)}x$
for all $x\in X\backslash \{w,z\}$. We have $g(q)=z$, and $g(p)\in \{w,z\}$
by Step 2. Therefore, $g(p)=z$ by strategy proofness. Create profile $s$
from $p$ by switching $w$ and $z$ in $p(2)$, with $s(i)=p(i)=q(i)$ for all $%
i\neq 2$. Step 2 implies $g(s)\in \{w,z\}$. If $g(s)=w$ then we have a
profile at which $w$ is at the top of person 1's ordering and $w$ is
selected. But then Lemma 6-2 and Step 3 of this lemma imply that $w$ is
selected whenever it is at the top of person 1's ordering, contradicting $%
g(q)=z$. Therefore, $g(s)=z$. Switch w and z in s(1) to create profile t for
which z is at the top of person 1's ordering and z is selcted (by
strategy-proofness).\textbf{\medskip } $\ \ \ \ \ \square $

\qquad Lemmas 6-2 and 6-3 imply that person 1 dictates $g$, and so we have
shown $SP(n,m+1)$. $\ \ \ \ \ \square $\medskip

\section{Conclusion\protect\medskip}

\qquad Combining the basis result, Theorem 5-1 and Theorem 6-1, we
have\medskip

\qquad \textbf{Theorem 7-1}. For all $m,n\geq 3$, $SP(n,m)$.\medskip

\qquad Important as this result is, we are also interested in rules that
have less than full range. \ Consider the rule $g$ defined as follows, where
we assume $X$ is large: $m>3$. \ Fix three alternatives, $\{a,b,c\}=S$ in $X$%
. \ At profile $u$, let $\psi (u)$ be the highest ranked alternative from $S$
according to ordering $u(1)$. \ Then $g(u)=\psi (u)$ unless there is a
unique $x\in X$ that Pareto-dominates $\psi (u)$, in which case $g(u)=x$. \
This is a UBM rule with full range. \ But the restriction $g^{\ast }$ of $g$
to $NP$ has range that is only $S$. \ This rule $g^{\ast }$ is then of less
than full range but is dictatorial on $NP$. \ That is a direct consequence
of Theorem 7-1 and the Equivalence Theorem from Section 3.\medskip

\qquad \textbf{Theorem 7-2.} Let $g^{\ast }$ be a strategy-proof rule on $%
NP(n,m)$ for $n\geq 3$ and $m\geq 3$ with Range($g^{\ast }$) = $S$ and $|S|$ 
$\geq 3$. Then $S$ is dictatorial.\medskip

\qquad \textbf{Proof:} For a fixed ordering $x,y,...,z$ of $X\backslash S$,
define a rule $g^{\ast \ast }$ on $NP(n,|S|)$ as follows: for profile $u$ in 
$NP(n,|S|)$, construct profile $u^{\ast }$ in $NP(n,m)$ as\medskip 
\[
u^{\ast }:%
\begin{tabular}{|c|c|c|c|}
\hline
\emph{1} & \emph{2} & \emph{\ldots } & \emph{n} \\ \hline
$x$ & $\vdots $ &  & $\vdots $ \\ 
$y$ & $u(2)$ & $\emph{\ldots }$ & $u(n)$ \\ 
$\vdots $ & $\vdots $ &  & $\vdots $ \\ 
$z$ &  &  &  \\ 
$\vdots $ & $z$ & $\emph{\ldots }$ & $z$ \\ 
$u(1)$ & $\vdots $ &  & $\vdots $ \\ 
$\vdots $ & $y$ &  & $y$ \\ 
& $x$ &  & $x$ \\ \hline
\end{tabular}%
\]%
\medskip Then define $g^{\ast \ast }(u)=g^{\ast }(u^{_{\ast }})$. It is easy
to show that $g^{\ast \ast }$ is strategy-proof. By the Equivalence Theorem, 
$g^{\ast \ast }$ is of full range: Range($g^{\ast \ast }$) $=S$. By Theorem
7-1, then, $g^{\ast \ast }$ is dictatorial. Then another application of the
Equivalence Theorem shows $g^{\ast }$ is dictatorial. \ \ $\square $\medskip

\section{References\protect\medskip}

Aswal, N, S Chatterji, and A Sen (2003). "Dictatorial domains," \textit{Econ
Theory} 22: 45-62.\medskip

Barber\`{a} S (1980) "Pivotal voters: A new proof of Arrow's theorem," 
\textit{Econ Letters} 6: 13-16.\medskip

Barber\`{a} S (1983a) "Pivotal voters: A simple proof of Arrow's theorem,"
In: Pattanaik PK and Salles M (eds) \textbf{Social Choice and Welfare}
Amsterdam; North-Holland, 31-35.\medskip

Barber\`{a} S (1983b) "Strategy-proofness and pivotal voters: A direct proof
of the Gibbard-Satterthwaite theorem," \textit{Int Econ Rev} 24:
413-418.\medskip

Barber\`{a} S (2001) "An introduction to strategy-proof social choice
functions," \textit{Soc Choice Welf} 18: 619-653.\medskip

Beno\^{\i}t JP (2000) "The Gibbard-Satterthwaite theorem: A simple proof," 
\textit{Econ Letters} 69: 319-322.\medskip

Campbell DE, Kelly JS (2014) "Universally beneficial manipulation: a
characterization," \textit{Soc. Choice Welf}. 43: 329-355.\medskip

Campbell DE, Kelly JS (2014b) "Two Theorems on the Range of Strategy-proof
Rules on a Restricted Domain.\textquotedblright\ \ 

\qquad arXiv: http://arxiv.org/abs/1408.7044\medskip

Gibbard A (1973) \textquotedblleft Manipulation of voting schemes: A general
result," \textit{Econometrica} 41: 587-601.\medskip

Larsson BS, Svensson LG (2006) "Strategy-proof voting on the full preference
domain," \textit{Math Soc Sci} 52: 272-287.\medskip\ 

Reny PJ (2001) "Arrow's theorem and the Gibbard-Satterthwaite theorem: A
unified approach," \textit{Econ Letters} 70: 99-105.\medskip

Satterthwaite M (1975) "Strategyproofness and Arrow's conditions: Existence
and correspondence theorems for voting procedures and social welfare
functions,"\ \textit{J Econ Theory} 10: 187-317.\medskip

Schmeidler D, Sonnenschein HF (1978)" Two proofs of the
Gibbard-Satterthwaite theorem on the possibility of a strategy-proof social
choice function," In: Gottinger HW and Leinfeller W (eds) \textbf{Decision
Theory and Social Ethics, Issues in Social Choice} Dordrecht; D. Reidel
227-234.\medskip

Sen A (2001) "Another direct proof of the Gibbard--Satterthwaite theorem," 
\textit{Econ Letters} 70: 381-385.\medskip

Svensson LG (1999) "The proof of the Gibbard--Satterthwaite theorem revisite%
\'{d}". University of Lund working paper:

\bigskip \bigskip 

Campbell: Department of Economics and The Program in Public Policy; 

The College of William and Mary; Williamsburg, VA 23187-8795, USA

E-mail: decamp@wm.edu\bigskip

Kelly: Department of Economics; Syracuse University; Syracuse, NY 

13244-1020, USA

E-mail: jskelly@maxwell.syr.edu

\end{document}